%% file: main.tex
\documentclass[a4paper,11pt, twoside]{article}

\usepackage{a4}
\usepackage{geometry} 
\usepackage{amsmath}
\usepackage{amssymb}
\usepackage{amsthm}
\usepackage{graphicx}
\usepackage{caption,subcaption}
\usepackage{multirow}
\usepackage{xstring}
\usepackage[utf8]{inputenc}
\usepackage{hyperref}
\usepackage{amstext,amsbsy,amsopn,eucal,enumerate}
\usepackage{tikz}
\usepackage{pgfplots}
\usepackage{tabularx}

\usepackage{soul}
\usepackage{algorithm} 
\usepackage{algpseudocode}

\def\R{\mathbb{R}}

\newcommand{\expn}{\operatorname{e}}

\newcommand{\diag}{\operatorname{diag}}

\newcommand{\kernel}{\operatorname{ker}}

\newcommand{\im}{\operatorname{im}}
\newcommand{\eig}{\operatorname{eig}}

\newcommand{\beq}{\begin{equation}}
\newcommand{\eeq}{\end{equation}}
\newcommand {\mat}      [1] {\left[\begin{array}{#1}}
\newcommand {\rix}          {\end{array}\right]}
\newcommand {\smat}      [1] {\left[\begin{smallmatrix}{#1}}
\newcommand {\srix}          {\end{smallmatrix}\right]}
\newcommand {\s}      [1] {\begin{smallmatrix}{#1}}
\newcommand {\se}          {\end{smallmatrix}}
\newcommand{\trace}{\operatorname{tr}}
\newcommand{\rank}{\operatorname{rank}}

\newcommand{\coreifive}{Core\texttrademark i5}

\newcommand{\cE}{\ensuremath{\mathcal{E}}}

\newtheorem{defn}{Definition}[section]
\newtheorem{remark}{Remark}

\newtheorem{example}[defn]{Example}
\newtheorem{lem}[defn]{Lemma}
\newtheorem{prop}[defn]{Proposition} 
\newtheorem{kor}[defn]{Corollary}
\newtheorem{thm}[defn]{Theorem}

\def\addlegendimage{\csname pgfplots@addlegendimage\endcsname}

\newlength\fheight
\newlength\fwidth

  \newcommand{\matlab}{MATLAB\textsuperscript{\textregistered}}
  \newcommand{\intel}{Intel\textsuperscript{\textregistered}}

\title{Full state approximation by Galerkin projection reduced order models for stochastic and bilinear systems}
\author{Martin Redmann\thanks{Martin Luther University Halle-Wittenberg, Institute of Mathematics, Theodor-Lieser-Str. 5, 06120 Halle (Saale), Germany, Email: {\tt 
martin.redmann@mathematik.uni-halle.de}.}\and Igor Pontes Duff \thanks{Max Planck Institute for Dynamics of Complex Technical Systems, Magdeburg, Germany, Email: {\tt 
pontes@mpi-magdeburg.mpg.de}.}
}

\begin{document}

\maketitle

\begin{abstract}
In this paper, the problem of full state approximation by model reduction is studied for stochastic and bilinear systems. 
Our proposed approach relies on identifying the dominant subspaces based on the reachability Gramian of a system. 
Once the desired subspace is computed, the reduced order model is then obtained by a Galerkin projection. 
We prove that, in the stochastic case, this approach either preserves mean square asymptotic stability or leads to reduced models whose minimal realization is mean square asymptotically stable. 
This stability preservation guarantees the existence of the reduced system reachability Gramian which is the basis for the full state error bounds that we derive. 
This error bound depends on the neglected eigenvalues of the reachability Gramian and hence shows that these values are a good indicator for the expected error in the dimension reduction procedure.
Subsequently, we establish the stability preservation result and the error bound for a full state approximation to bilinear systems in a similar manner. These latter results are based on a recently proved link between
stochastic and bilinear systems. We conclude the paper by numerical experiments using a benchmark problem. We compare this approach with balanced truncation and show that it performs well in reproducing the full state of the system.
\end{abstract}



\section{Introduction}

Galerkin approximation is an important methodology to obtain surrogate models for high fidelity systems. It relies on the fact that, in many applications, the state of the system is well approximated in a
lower-dimensional subspace. In other words, for $x(t) \in \R^n$ with $n \gg 1$, there exist $V \in \R^{n\times r}$ such that $x(t) \approx V\hat{x}(t)$. 
The choice of the right basis $V$ plays a crucial role in the approximation quality. Several approaches to construct the dominant subspaces have been proposed for deterministic linear systems, see, e.g., \cite{benner2015survey}. 
Among them, a commonly used approach is the proper orthogonal decomposition (POD) \cite{berkooz1993proper,kunisch2001galerkin}, which identifies dominant subspaces empirically by extracting
them from snapshot matrices. These empirical methods have been successfully used in many applications. 
However, they are \emph{input-dependent} in the setup of  control systems, i.e., the quality of reduced order model (ROM) will depend on the choices of inputs used to generate the snapshots. 
In the setup of stochastic systems considered in this paper, a POD approach would require the simulation of an enormous amount of samples being numerically costly in practice. Moreover, a deeper theoretical analysis 
of such empirical methods is often not feasible. \smallskip

Given stability in the original model, it is of interest to preserve this property in the reduced system. 
Galerkin methods have been shown to preserve stability for deterministic linear dissipative systems, see, e.g., \cite{selga2012stability}. 
Additionally, the authors in  \cite{prajna2003pod} propose a POD scheme combined with linear matrix inequalities to construct ROMs that are locally stable in a nonlinear deterministic setting.
In general, Krylov based methods for deterministic linear \cite{freund2003model} and bilinear systems \cite{breiten2010krylov} do not guarantee stability in the ROM. 
However, in \cite{morPul19}, the authors have proposed equivalent dissipative realizations to arbitrary Galerkin projected linear systems that are stable. 
To the authors' knowledge, stability preservation using Galerkin projections has not been studied in the literature of stochastic systems.\smallskip

In this work, we focus on model order reduction of linear stochastic and bilinear systems aiming for a full state approximation. 
Therefore, we study a Galerkin approach based on the dominant reachability subspaces, which leads to input-independent projections, i.e., 
the corresponding ROMs are (pathwise) accurate for a large set of inputs. This approach relies on the computation of the reachability Gramian of the underlying dynamical system.
These Gramians are encoded by generalized Lyapunov equations and, hence, can be computed in a numerically efficient way in large-scale settings (see, e.g., the review papers \cite{BenS13, Sim16a} for low-rank methods).
Once the Lyapunov equation is solved, the dominant subspaces are then identified by the span of eigenvectors of the Gramian associated with the large eigenvalues. 
Hence, the ROM is obtained by projecting the dynamical system onto the identified subspace. 
We show that this procedure either preserves the underlying stability or at least ensures the existence of the reduced order Gramian, which is vital for the error analysis. As a consequence, the minimal realization of the reduced model is stable.
Subsequently, we propose error bounds for the approximation, which show how the reduction error is related to the neglected eigenvalues of the Gramian. 
It is worth noticing that this approach has already been successfully applied in the literature of deterministic linear time-invariant systems, e.g., in the context of structured systems \cite{sorensen2005model},
and port-Hamiltonian systems \cite{polyuga2010model}.
However, to the authors’ knowledge, the stability analysis and error bounds for stochastic and bilinear systems considered in this paper have not even been established for deterministic linear systems so far.\smallskip

It is worth mentioning that one way to address the problem of full state approximation is to use balanced truncation, see \cite{moo1981} for deterministic linear systems and 
\cite{typeIBT, bennerdamm, redmannbenner, redmannPhD} for stochastic and bilinear systems. This method is generally suitable when one wants to approximate a quantity of interest $y(t) = Cx(t)$, with 
$C\in \R^{p \times n}$ and $p \ll n$. 
For this method, one needs to compute the observability Gramian in addition to the reachability Gramian. Subsequently, a ROM is obtained by Petrov-Galerkin projection based on these two Gramians. 
One advantage of balanced truncation is that it is, under some mild conditions, stability preserving \cite{redbendamm, pernebo1982model} and it guarantees error bounds \cite{redmannbenner, enns1984model, h2_bil}. 
However, whenever,  $p \approx n$, it suffers from the issue that the computation 
of the observability Gramian is not feasible in practice since low-rank methods are no longer applicable in this context. This scenario is given if the full state shall be approximated, since $C =I$ in this case. 
For the deterministic linear case, the authors in \cite{morBenS10} propose a scheme enabling the computation of the Petrov-Galerkin projection without the explicit computation of the observability Gramian. 
The approach is numerically feasible but costly since it relies on a quadrature scheme using the low-rank factors of the reachability Gramian pre-multiplied by shifted systems. Additionally, those results are no longer 
applicable for stochastic and bilinear systems. \smallskip

The paper is organized as follows. In Section \ref{sec:StochBackground}, we present the main setup for linear stochastic systems and the concept of mean square asymptotic stability along with some literature results. 
Then, in Section \ref{romsec}, we described the proposed Galerkin projection based procedure using the reachability Gramian. Additionally, an interpretation of the dominant subspaces is derived therein. 
Section \ref{sec:PropROM} is dedicated to showing the properties of the ROM.  First, we prove that this procedure either constructs a reduced model which is mean square asymptotically stable 
or a ROM which has a realization satisfying the desired stability property. It is worth noticing that modified versions of those results are also valid for
the class of deterministic bilinear systems, which we also establish in this paper. In Subsection \ref{subsec:ErrBound}, we derive bounds for the approximation error and their relation to the neglected singular values of the 
reachability Gramian. In Section \ref{sec:MORbil}, similar results for the class of bilinear systems, including stability preservation and error bounds, are derived. 
In Section \ref{sec:NumExp}, some numerical experiments are conducted to illustrate the performance of the proposed approach and in order to compare it with balanced truncation.

\section{Linear stochastic systems and mean square stability}\label{sec:StochBackground}
\subsection{Stochastic problem setup}
We consider the following linear stochastic systems 
\begin{align}\label{stochstate}
 dx(t) = [Ax(t) +  B u(t)]dt + \sum_{i=1}^{q}  N_i x(t) dW_i(t),\quad t\geq 0,
\end{align}
where we assume that $A, N_i \in \R^{n\times n}$ and $B \in \R^{n\times m}$ are constant matrices and its initial condition is $x(0)=0$. 
The vectors $x$ and $u$ are called state and control input, respectively. Moreover, let 
$W=\left(W_1, \ldots, W_{q}\right)^\top$ be an $\mathbb R^{q}$-valued standard Wiener process for simplicity of the notation. The results can be extended to square integrable L\'evy processes with mean zero and
general covariance matrix (see, e.g., \cite{redmannPhD}). All stochastic processes appearing in this paper, are defined on a filtered probability space
$\left(\Omega, \mathcal F, (\mathcal F_t)_{t\geq 0}, \mathbb P\right)$\footnote{$(\mathcal F_t)_{t\geq 0}$ is right continuous and complete.}. In addition, $W$ is $(\mathcal F_t)_{t\geq 0}$-adapted and its 
increments $W(t+h)-W(t)$ are independent of $\mathcal F_t$ for $t, h\geq 0$. Throughout this paper, we assume that $u$ is an $(\mathcal F_t)_{t\geq 0}$-adapted control that is square integrable, meaning that 
\begin{align*}
\left\| u\right\|^2_{L^2_T} := \mathbb E\int_0^T \left\| u(s)\right\|_2^2 ds < \infty
\end{align*}
for all $T>0$, where $\left\| \cdot\right\|_2$ denotes the Euclidean norm. Moreover,  $\left\| \cdot\right\|_F$ will denote the Frobenius norm, whereas $\left\| \cdot\right\|$ represents an arbitrary matrix/vector norm.  
The aim is to identify a low-dimensional subspace $\mathcal V$ of $\mathbb R^n$ that approximates the manifold 
of the state $x$. Choosing a matrix $V\in\mathbb R^{n\times r}$ of orthonormal basis vectors of $\mathcal V$, an approximation of the form $V \hat x(t) \approx x(t)$ can be constructed. Inserting this approximation into the original system \eqref{stochstate}, we enforce a Petrov-Galerkin condition 
by multiplying the residual with $V^\top$ leading to a ROM 
\begin{align}\label{red_stochstate_mul}
 d\hat x(t) = [\hat A\hat x(t) +  {\hat B} u(t)]dt + \sum_{i=1}^{q} {\hat N}_i \hat x(t) dW_i(t), \quad t\geq 0,
\end{align}
where $\hat A = V^\top A V$, $\hat B = V^\top B$, $\hat N_i = V^\top\hat N_i V$ and $\hat x(t)\in\mathbb R^r, \hat x(0)=0$, with $r\ll n$. 
Our main goal is to construct the matrix $V$, such that, the approximation error is small for every input $u$ considered.

\subsection{Mean square asymptotic stability and generalized Lyapunov operators}\label{sec_stab}

We introduce the fundamental solution $\Phi$ to \eqref{stochstate}. It is defined as the $\R^{n\times n}$-valued solution to 
\begin{align}\label{funddef}
 \Phi(t, s)=I+\int_s^t A \Phi(\tau, s) d\tau+\sum_{i=1}^{q} \int_s^t N_i \Phi(\tau, s)dW_i(\tau), \quad t\geq s.
\end{align}
It is the operator that maps the initial condition $x_0$ to the solution of the homogeneous state equation, i.e., $u \equiv 0$, with initial time $s\geq 0$. We additionally define $\Phi(t) := \Phi(t, 0)$. Moreover, notice that we have 
$\Phi(t, s) = \Phi(t)\Phi^{-1}(s)$.\smallskip

Throughout this paper, we assume that the uncontrolled state equation \eqref{stochstate} is mean square asymptotically stable, i.e, $\mathbb E\left\|\Phi(t)\right\|^2 \lesssim \expn^{-c t}$
for some constant $c>0$. With $\lambda(\cdot)$ denoting the spectrum of a matrix/operator, this is equivalent to 
 \begin{align}\label{assumptionstab}
  \lambda\left(K\right)\subset \mathbb C_-, \end{align}
where $K:=I\otimes A + A\otimes I+ \sum_{i=1}^{q} N_i \otimes N_i$ and $\cdot\otimes\cdot$ is the Kronecker product of two matrices, see for instance \cite{damm, staboriginal}. 
Moreover, the system is called mean square stable if $  \lambda\left(K\right)\subset \overline{\mathbb C_-}$. 
Notice that $\lambda(K)=\lambda(\mathcal L_A + \Pi_N)$, where the generalized Lyapunov operator $\mathcal L_A + \Pi_N$ is defined by
$X\mapsto \mathcal L_A(X) = AX + X A^\top$ and $X\mapsto \Pi_N(X) = \sum_{i=1}^q N_i X N_i^\top$. In Section \ref{romsec}, 
we will introduce a reduced system which does not necessarily preserve \eqref{assumptionstab} but it is always mean square stable, i.e., $K$ 
can additionally have eigenvalues on the imaginary axis. Therefore, we need the following sufficient conditions for mean square stability.
\begin{lem}\label{lem_stab_plus_zero}
Given a matrix $Y\geq 0$, let us assume that there exists $X>0$ such that \begin{align*}
            \mathcal L_A(X) + \Pi_N(X)\leq -Y.                                                               
                                                                          \end{align*}
Then, we have $\lambda\left(K\right)\subset \overline{\mathbb C_-}$.
\end{lem}
\begin{proof}
An algebraic proof can be found in \cite[Corollary 3.2]{redbendamm}. We refer to \cite[Lemma 6.12]{redmannPhD} for a probabilistic approach.
\end{proof}
Let $\alpha(\mathcal L_A + \Pi_N):=\max\{\Re(\mu): \mu\in\lambda(\mathcal L_A + \Pi_N)\}$ be the spectral abscissa of the operator $\mathcal L_A + \Pi_N$, with $\Re(\cdot)$ being the real part of a complex number.
Since the stability of a stochastic system is related to the eigenvalues of $\mathcal L_A + \Pi_N$, we formulate the following result.
\begin{lem}\label{lem_zero_eig}
 There exists $V_1\geq 0$, $V_1\neq 0$, such that $\mathcal L_A(V_1) + \Pi_N(V_1)= \alpha(\mathcal L_A + \Pi_N) V_1$.
\end{lem}
\begin{proof}
 A proof in a more general framework can be found in \cite[Section 3.2]{damm}. We also refer to \cite[Theorem 3.1]{redbendamm} and the references therein.
\end{proof}
Finally, spectral properties of the Kronecker matrix involving both the reduced and the original model matrices are required. 
\begin{lem}\label{mixed_stab}
  Given that the full model is mean square asymptotically stable, whereas the reduced system is just mean square stable, i.e., we have    \begin{align*}
          \lambda\left(I\otimes A + A\otimes I+ \sum_{i=1}^{q} N_i \otimes N_i\right)\subset \mathbb C_-\text{ and }  
          \lambda\left(I\otimes \hat A + \hat A\otimes I+ \sum_{i=1}^{q} \hat N_i \otimes \hat N_i\right)\subset \overline{\mathbb C_-}.
                           \end{align*}
Then, it holds that      \begin{align*}
          \lambda\left(I\otimes A + \hat A\otimes I+ \sum_{i=1}^{q} \hat N_i \otimes N_i\right)\subset \mathbb C_-.
                           \end{align*}                                                                                                              
                                                                                                                           \end{lem}
    \begin{proof}
A probabilistic version can be found in \cite[Lemma 6.12]{redmannPhD} and an algebraic approach is given in \cite[Proposition 3.4]{redbendamm}.
    \end{proof}

\section{Dominant subspaces and reduced order model}\label{romsec}

\subsection{Dominant subspaces of \eqref{stochstate}}\label{domsubsec}

We identify the redundant information in the system by using the reachability Gramian \[P:= \mathbb E \int_0^\infty \Phi(s)BB^\top \Phi^\top (s) ds.\] Notice that $P$ exists due to the exponential decay of the fundamental solution $\Phi$. 
Practically, one can compute $P$ by solving a generalized Lyapunov equation. Using Lemma \ref{lemdgl} with $s=0$, $\hat A = A$, $\hat B = B$, $\hat N_i = N_i$ and $t\rightarrow \infty$, we obtain that the reachability Gramian is a solution of the following generalized Lyapunov equation
\begin{align}\label{reach_gram}
 A P + P A^\top + \sum_{i=1}^q N_i P N_i^\top = -B B^\top.
\end{align}
Let $x(t, x_0, u)$, $t\geq 0$, denote the solution of \eqref{stochstate} with initial value $x_0$ and control $u$. Then, for $z\in\mathbb R^n$, we have
\begin{align}\label{interpgram1}
\sup_{t\in[0, T]}\mathbb E \left\vert\langle x(t, 0, u), z \rangle_2\right\vert \leq 
\left(z^\top P z\right)^{\frac{1}{2}} \left\|u\right\|_{L^2_T}
\end{align}
using the results in \cite{redmannspa2}. Let $(p_k)_{k=1,\ldots, n}$ be an orthonormal basis of $\mathbb R^n$ consisting of eigenvectors of $P$. Then, the state variable can be written as \begin{align*}
 x(t, 0, u)=\sum_{k=1}^n  \left\langle x(t, 0, u), p_{k} \right\rangle_2 p_{k}.  \end{align*}
Setting $z=p_k$ in \eqref{interpgram1}, we obtain
 \begin{align}\label{diffreachjanein}
\sup_{t\in[0, T]}\mathbb E \left\vert\langle x(t, 0, u), p_{k}  \rangle_2\right\vert \leq \lambda_{k}^{\frac{1}{2}}\left\|u\right\|_{L^2_T},
\end{align}
where $\lambda_{k}$ is the corresponding eigenvalue. Consequently, we see that the direction $p_k$ is completely irrelevant if $\lambda_{k} = 0$. On the other hand, if $\lambda_{k}$ is not zero but small, then 
a large component in the direction of $p_{k}$ requires a large amount of energy by \eqref{diffreachjanein}. Therefore, the eigenspaces of $P$ belonging to the small eigenvalues can also be neglected.\smallskip

A ROM can now be obtained by removing the unimportant subspaces from \eqref{stochstate}. This is done by first diagonalizing $P$. 
If $P$ is diagonal, we have that $p_k$ is the $k$th unit vector 
and the diagonal entries of $P$ indicate the relevance of the respective unit vector. A reduced system can then be easily derived by truncating the components of $x$ associated to the small/zero entries $\lambda_{k}$ 
of a Gramian of the form $P= \diag(\lambda_1, \ldots, \lambda_n)$.

\subsection{Reduced order model by Galerkin projection}\label{rom_computattion}

We introduce the eigenvalue decomposition of the reachability Gramian as follows\begin{align*}
          P= S^\top \Lambda S,                                    
                                             \end{align*}
where $S^{-1}=S^\top$ and $\Lambda= \smat{\Lambda}_{1}& 0\\ 
0 &{\Lambda}_{2}\srix =  \diag(\lambda_1, \ldots, \lambda_n)$ is the matrix of eigenvalues of $P$. For simplicity, let us assume that the spectrum of $P$ is ordered, i.e., $\lambda_1\geq \ldots \geq \lambda_n\geq 0$ so that 
$\Lambda_2$ contains the small eigenvalues.  
Let us do a state space transformation using the matrix $S$. The transformed state variable then is $x_b = S x$. Plugging this into \eqref{stochstate}, we find 
\begin{equation}\label{balancedsystem}
\begin{aligned}
 dx_b(t) &= [A_bx(t) +  B_b u(t)]dt + \sum_{i=1}^{q}  N_{i, b} x(t) dW_i(t),\quad t\geq 0,\\
 x(t) &= S^\top x_b(t),
 \end{aligned}
\end{equation} 
where the balanced matrices are given by
\begin{align}\label{partition}
A_b:=S{A}S^\top= \smat{A}_{11}&{A}_{12}\\ 
{A}_{21}&{A}_{22}\srix,\quad B_b:= S{B} = \smat{B}_1\\ {B}_2\srix,\quad N_{i, b} := S{N_i}S^\top= \smat{N}_{i, 11}&{N}_{i, 12}\\ 
{N}_{i, 21}&{N}_{i, 22}\srix.\end{align}
We refer to \eqref{partition} as the balanced realization of the linear stochastic system.
The fundamental solution of the balanced realization is $\Phi_b = S \Phi S^\top$ which can be seen by multiplying \eqref{funddef} with $S$ from the left and with $S^\top$ from the right. Therefore, the reachability Gramian of 
\eqref{balancedsystem} is \begin{align}\label{balanced_gram}
   P_b:= \mathbb E \int_0^\infty \Phi_b(s)B_bB_b^\top \Phi_b^\top (s) ds = S P S^\top = \Lambda.
                          \end{align}
We partition $x_b=\smat x_1\\ x_2 \srix$, where $x_1$ and $x_2$ are associated to $\Lambda_1$ and $\Lambda_2$, respectively. Now, exploiting the insights of Section \ref{domsubsec}, $x_2$ barely contributes to the system 
dynamics. We obtain the reduced system by truncating the  equation related to $x_2$ in \eqref{balancedsystem}. Furthermore, we set the remaining $x_2$ components equal to zero. This yields a reduced 
system \eqref{red_stochstate_mul} with matrices \begin{align}\label{red_sys_one_sided_BT}
\hat A = A_{11} = V^\top A V, \quad \hat B = B_1 = V^\top B,\quad   \hat N_{i} =   {N}_{i, 11} = V^\top N_i V,                                                                   
                                                                            \end{align}
where $V$ are the first $r$ columns of $S^\top =  \smat V & S_{2}\srix$.  \smallskip

In large-scale settings, the reachability Gramian can be computed using low-rank methods (see \cite{BenS13, Sim16a}), i.e., 
we find a matrix  $Z_P\in \R^{n \times l}$, with $l \ll n$, such that $P \approx Z_PZ_P^\top $. Consequently, in this setup, the Galerkin projection can be identified using the singular value decomposition
of  $Z_P$.

\section{Properties of the reduced system}\label{sec:PropROM}

\subsection{Mean square stability and reduced order Gramian}

In this section, we study stability preservation and the existence of the Gramian for the reduced system in \eqref{red_sys_one_sided_BT}. The next result guarantees mean square stability.
\begin{prop}\label{mean_square_stable_red_sys}
Suppose that $\Lambda_1 = \diag(\lambda_1, \ldots, \lambda_r)>0$. Then, the reduced order system \eqref{red_stochstate_mul} with $\hat A = A_{11}$ and $\hat N_{i} =   {N}_{i, 11}$, introduced in \eqref{partition}, is mean square stable, i.e., \begin{align*}
\lambda(I\otimes A_{11} + A_{11}\otimes I+ \sum_{i=1}^{q} N_{i, 11} \otimes N_{i, 11})\subset \overline{\mathbb C_-}.                                                                                                                                                                          
                                                                                                                                                                         \end{align*}
\end{prop}
\begin{proof}
 According to \eqref{balanced_gram}, the balanced reachability Gramian is the diagonal matrix $\Lambda$ of eigenvalues of $P$. Using the partition of the balanced matrices in \eqref{partition}, we therefore have 
 \begin{align*}
   &\smat{A}_{11}&{A}_{12}\\ {A}_{21}&{A}_{22}\srix \smat{\Lambda}_{1}& 0\\ 
0 &{\Lambda}_{2}\srix +\smat{\Lambda}_{1}& 0\\ 
0 &{\Lambda}_{2}\srix \smat{A}_{11}^\top&{A}_{21}^\top\\ {A}_{12}^\top&{A}_{22}^\top\srix+
\sum_{i=1}^q \smat{N}_{i, 11}&{N}_{i, 12}\\ {N}_{i, 21}&{N}_{i, 22} \srix \smat{\Lambda}_{1}& 0\\ 
0 &{\Lambda}_{2}\srix  \smat{N}_{i, 11}^\top &{N}_{i, 21}^\top\\ {N}_{i, 12}^\top&{N}_{i, 22}^\top \srix \\
   &=-\smat{B}_1B_1^\top&{B}_1B_2^\top\\ {B}_2 B_1^\top & {B}_2B_2^\top\srix. 
 \end{align*}
The left upper block of the above equation is  \begin{align*}
             A_{11} \Lambda_1 + \Lambda_1 A_{11}^\top + \sum_{i=1}^q N_{i, 11} \Lambda_1 N_{i, 11}^\top = -B_1B_1^\top - \sum_{i=1}^q N_{i, 12} \Lambda_2 N_{i, 12}^\top \leq 0.
                  \end{align*}
Since $\Lambda_1>0$ by assumption, Lemma \ref{lem_stab_plus_zero} yields the claim.
\end{proof}
Using Proposition \ref{mean_square_stable_red_sys} and Lemma \ref{lem_zero_eig}, the reduced order system is asymptotically mean square stable if and only if 
$0\not\in \lambda(I\otimes A_{11} + A_{11}\otimes I+ \sum_{i=1}^{q} N_{i, 11} \otimes N_{i, 11})$. With the following example it is shown that the zero eigenvalue can indeed occur.
\begin{example}\label{example1}
  Let $N_i = 0$, $A = \smat   0& -10 \\ 1 & -10 \srix$ and $B = \smat 0\\10\srix$. Then, system \eqref{stochstate} is already balanced 
since the reachability Gramian is given by
\[ 
P = \int_0^{\infty} \expn^{A s}BB^\top \expn^{A^\top s}ds = \smat 50 & 0 \\ 0 & 5 \srix.
\]
Moreover, we have $\rank(\smat B & A B\srix) = 2$ which means that the system is reachable or locally reachable if $N_i$ were non zero. According to Section \ref{rom_computattion} the reduced matrices are 
$A_{11} = 0$, $B_1 = 0$ and $N_{i, 11} = 0$. Consequently, the reduced system is not asymptotically stable and the reachability of the system is also lost since the reduced system is uncontrolled.
\end{example}
Example \ref{example1} shows a difference to balanced truncation for stochastic systems, where mean square asymptotic stability is preserved under relatively general conditions \cite{redbendamm, bennerdammcruz}. 
Mean square asymptotic stability ensures the existence of the reachability Gramian $\hat P:= \mathbb E \int_0^\infty \hat \Phi(s)B_1B_1^\top \hat \Phi^\top (s) ds$, where $\hat \Phi$ represents the fundamental solution of the 
reduced system. However, asymptotic stability is only a sufficient condition for $\hat P$ to exist. The next example illustrates such a scenario. 
\begin{example}\label{example2}
 Let $N_i = 0$, $A=\smat -1& -1\\ -1 &-1\srix$ and $B= \smat 1\\ 1\srix$. Then, we have \begin{align*}
   \Phi(t) = \expn^{At} =   \smat  {\frac{\expn^{-2 t}+1}{2}}  & {\frac{\expn^{-2 t}-1}{2}}  \\ {\frac{\expn^{-2 t}-1}{2}}  & {\frac{\expn^{-2 t}+1}{2}}  \srix   \quad \text{and}\quad \Phi(t)B =  \smat  \expn^{-2 t}   \\ \expn^{-2 t} \srix.                                                                  
                                                                                         \end{align*}
Clearly, $\Phi$ does not decay exponential to zero but $\Phi B$ does. Therefore, the reachability Gramian exist and is $P= \smat 0.25 & 0.25\\ 0.25 & 0.25\srix$, whereas the set of solutions to \eqref{reach_gram} is given by 
$P + y \smat -1 & 1\\ 1 & -1\srix$, $y\in\mathbb R$.
\end{example}
Example \ref{example2} emphasizes that it is important to distinguish between a Gramian (given by an integral representation) and a solution of a Lyapunov equation. 
The next theorem proves that even if the reduced system is not mean square asymptotically stable, the existence of the reduced order reachability Gramian can be guaranteed. 
This is one of the main results of the paper which is also vital for later considerations, where we prove an error bound in which $\hat P$ is involved.
\begin{thm}\label{thm_red_gram_exisits}
 Given the reduced system \eqref{red_stochstate_mul} with matrices $\hat A = A_{11}$, $\hat B = B_1$, $\hat N_{i} =   {N}_{i, 11}$ defined in \eqref{partition} and $\Lambda_1 = \diag(\lambda_1, \ldots, \lambda_r)>0$. 
 Moreover, let $\hat \Phi$ denote the fundamental 
 solution to this reduced order system. Then, there is a constant $c>0$ such that $\mathbb E\left\|\hat \Phi(t) B_1\right\|_F^2 \lesssim \expn^{-c t}$. Hence, the reachability Gramian 
 $\hat P:= \mathbb E \int_0^\infty \hat \Phi(s)B_1B_1^\top \hat\Phi^\top (s) ds$ exists and satisfies 
                                         \begin{align}\label{lyap_red_sys}
 A_{11} \hat P + \hat P A_{11}^\top + \sum_{i=1}^q N_{i, 11} \hat P N_{i, 11}^\top = -B_1 B_1^\top.
\end{align}
\end{thm}
\begin{proof}
We set $\hat K:=I\otimes A_{11} + A_{11}\otimes I+ \sum_{i=1}^{q} N_{i, 11} \otimes N_{i, 11}$ and consider the case that the reduced system is mean square asymptotically stable, i.e., $0\not\in\lambda(\hat K)$. 
According to Section \ref{sec_stab} this is equivalent to $\mathbb E\left\|\hat \Phi(t)\right\|_F^2 \lesssim \expn^{-c t}$ implying $\mathbb E\left\|\hat \Phi(t) B_1\right\|_F^2 \lesssim \expn^{-c t}$. Given this condition 
 the infinite integral $\hat P$ exists. Moreover, using Lemma \ref{lemdgl} with $A=\hat A=A_{11}$, $B=\hat B=B_1$ and $N_{i}=\hat N_{i}=N_{i, 11}$ and exploiting  that the left hand side of \eqref{dglmixed} 
tends to zero if $t\rightarrow \infty$, we see that $\hat P$ solves \eqref{lyap_red_sys}.\smallskip

Let us consider the case of $0\in\lambda(\hat K) = \lambda({\hat K}^\top)$. If further $B_1=0$, the result of this theorem is true. Therefore, we additionally assume that $B_1 \neq 0$.
Then, by Lemma \ref{lem_zero_eig}, there exists $\hat V\geq 0$ such that \begin{align}\label{eig_vel_eq}
             \mathcal L_{A_{11}^\top}(\hat V) + \Pi_{N_{11}^\top}(\hat V) =  A_{11}^\top \hat V + \hat V A_{11} + \sum_{i=1}^q N_{i, 11}^\top \hat V N_{i, 11} = 0.                                                                                                                                
                                                                                                                                                  \end{align}
Moreover, according to the proof of Proposition \ref{mean_square_stable_red_sys}, we have \begin{align}\label{left_uper_block}
             A_{11} \Lambda_1 + \Lambda_1 A_{11}^\top + \sum_{i=1}^q N_{i, 11} \Lambda_1 N_{i, 11}^\top = -B_1B_1^\top - \sum_{i=1}^q N_{i, 12} \Lambda_2 N_{i, 12}^\top =:-R.
                  \end{align}
We observe that \begin{align*}
-\langle R, \hat V  \rangle_F = \langle \mathcal L_{A_{11}}(\Lambda_1) + \Pi_{N_{11}}(\Lambda_1), \hat V  \rangle_F = \langle \Lambda_1, \mathcal L_{A_{11}^\top}(\hat V) + \Pi_{N_{11}^\top}(\hat V)  \rangle_F = 0.
                \end{align*}
  Using the properties of the trace this yields \begin{align*}
 \left\|\hat V^{\frac{1}{2}} B_1\right\|_F^2 + \sum_{i=1}^q \left\|\hat V^{\frac{1}{2}} N_{i, 12} \Lambda_2^{\frac{1}{2}}\right\|_F^2   
 &= \trace(\hat V^{\frac{1}{2}} B_1 B_1^\top \hat V^{\frac{1}{2}}) +   \sum_{i=1}^q   \trace(\hat V^{\frac{1}{2}} N_{i, 12}\Lambda_2 N_{i, 12}^\top \hat V^{\frac{1}{2}})\\ &= \langle R, \hat V  \rangle_F = 0.
              \end{align*}                                                                                                                              
This implies that \begin{align}\label{zero_coef}
   \hat V B_1 = 0  \quad\text{and}\quad   \hat V N_{i, 12} \Lambda_2^{\frac{1}{2}}= 0.          
                  \end{align}
The case $\hat V>0$ is excluded since then it holds that $B_1 = 0$. Therefore, we consider the scenario in which $\hat V$ does not have full rank. We then assume that $\hat V$ is an
eigenvector with maximal rank, i.e., for any other eigenvector $\tilde V\geq 0$ corresponding to the zero eigenvalue, we have $\rank(\tilde V)\leq \rank(\hat V)$.\smallskip

Introducing the eigenvalue decomposition of $\hat V$:\begin{align}\label{eigdecVhat}
             \hat V = \smat {\hat V_1}&{\hat V_2} \srix  \smat{\hat D}&0\\0&0  \srix   \smat {\hat V_1^\top} \\ {\hat V_2^\top} \srix     = \hat V_1 \hat D \hat V_1^\top,                                                                                                                        
                                                                                                                                                \end{align}
$\hat D>0$, we find a basis of the kernel by the columns of $\hat V_2$, i.e., $\kernel(\hat V)= \im(\hat V_2)$. Inserting \eqref{eigdecVhat} into \eqref{zero_coef} yields
\begin{align}\label{zero_coef2}
   \hat V_1^\top B_1 = 0  \quad\text{and}\quad   \hat V_1^\top N_{i, 12} \Lambda_2^{\frac{1}{2}}= 0.          
                  \end{align}
We use a state space transformation based on $\hat S= \smat {\hat V_1^\top} \\ {\hat V_2^\top} \srix $ 
involving the following matrices\begin{equation}\begin{aligned}\label{part_trans_red_sys}
&\hat S A_{11} \hat S^\top=:\smat {\hat A_{11}}& {\hat A_{12}}\\ {\hat A_{21}}& {\hat A_{22}}\srix,\quad \hat S N_{i, 11} \hat S^\top =: \smat {\hat N_{i, 11}}& {\hat N_{i, 12}}\\ {\hat N_{i, 21}}& {\hat N_{i, 22}}\srix, \quad 
\hat S \Lambda_{1} \hat S^\top =: \smat {\hat P_{11}}& {\hat P_{12}}\\ {\hat P_{12}^\top}& {\hat P_{22}}\srix,\\
& \hat S B_1 =\smat {\hat V_1^\top B_1}\\{\hat V_2^\top B_1}  \srix = \smat 0 \\ {\hat V_2^\top B_1}  \srix, \quad 
\hat S N_{i, 12} \Lambda_2^{\frac{1}{2}} =\smat {\hat V_1^\top N_{i, 12} \Lambda_2^{\frac{1}{2}}}\\{\hat V_2^\top N_{i, 12} \Lambda_2^{\frac{1}{2}}}  \srix = \smat 0 \\ {\hat V_2^\top N_{i, 12} \Lambda_2^{\frac{1}{2}}}  \srix,
\end{aligned}\end{equation}
where \eqref{zero_coef2} was exploited. We multiply \eqref{left_uper_block} with $\hat S$ from the left and with $\hat S^\top$ from the right and obtain
\begin{align}\nonumber
 &\smat {\hat A_{11}}& {\hat A_{12}}\\ {\hat A_{21}}& {\hat A_{22}}\srix \smat {\hat P_{11}}& {\hat P_{12}}\\ {\hat P_{12}^\top}& {\hat P_{22}}\srix +
 \smat {\hat P_{11}}& {\hat P_{12}}\\ {\hat P_{12}^\top}& {\hat P_{22}}\srix \smat {\hat A_{11}^\top}& {\hat A_{21}^\top}\\ {\hat A_{12}^\top}& {\hat A_{22}^\top}\srix+
 \sum_{i=1}^q   \smat {\hat N_{i, 11}}& {\hat N_{i, 12}}\\ {\hat N_{i, 21}}& {\hat N_{i, 22}}\srix  \smat {\hat P_{11}}& {\hat P_{12}}\\ {\hat P_{12}^\top}& {\hat P_{22}}\srix 
  \smat {\hat N_{i, 11}^\top}& {\hat N_{i, 21}^\top}\\ {\hat N_{i, 12}^\top}& {\hat N_{i, 22}^\top}\srix\\ \label{parttransformedsys}
  &= -\smat 0& 0\\ 0& {\hat R}\srix 
\end{align}
with $\hat R = \hat V_2^\top B_1B_1^\top \hat V_2 + \sum_{i=1}^q  \hat V_2^\top N_{i, 12} \Lambda_2 N_{i, 12}^\top \hat V_2\geq 0$. Before we evaluate the blocks of \eqref{parttransformedsys}, we show that 
\begin{align} \label{rightupperzero}
 \hat A_{12} = \hat V_1^\top A_{11} \hat V_2 = 0\quad \text{and}\quad \hat N_{i, 12} = \hat V_1^\top N_{i, 11} \hat V_2 = 0.
\end{align}
To do so, we show that the kernel of $\hat V$ is invariant under multiplication with $A_{11}$ and $N_{i, 11}$. Let $z\in\kernel(\hat V)$. Then, we obtain \begin{align*}
   0 = z^\top\left(A_{11}^\top \hat V + \hat V A_{11} + \sum_{i=1}^q N_{i, 11}^\top \hat V N_{i, 11}\right) z   =     \sum_{i=1}^q z^\top N_{i, 11}^\top \hat V N_{i, 11} z
   = \sum_{i=1}^q \left\| \hat V^{\frac{1}{2}} N_{i, 11} z \right\|_2^2                                                                                                                                         
                                                                                                                                                          \end{align*}
implying that $\hat V N_{i, 11} z = 0$. Using this fact provides that  \begin{align*}
   0 = \left(A_{11}^\top \hat V + \hat V A_{11} + \sum_{i=1}^q N_{i, 11}^\top \hat V N_{i, 11}\right) z   =  \hat V A_{11} z.                                                                                                                                   
                                                                                                                                                          \end{align*}
Hence, we have $A_{11} \kernel(\hat V), N_{i, 11} \kernel(\hat V)\subset \kernel(\hat V)$. Since the columns of $\hat V_2$ span $\kernel(\hat V)$ and due to the invariance, there exist suitable matrices 
$\tilde A_{11}$ and $\tilde N_{i, 11}$ such that\begin{align}\label{conseq_invariance}
                                           A_{11} \hat V_2 = \hat V_2 \tilde A_{11}\quad\text{and}\quad N_{i, 11} \hat V_2 = \hat V_2 \tilde N_{i, 11}.
                                          \end{align}
Exploiting that $\hat V_1^\top \hat V_2 = 0$ gives us \eqref{rightupperzero}. Moreover, we see that $\hat A_{22} = \hat V_2^\top A_{11} \hat V_2 = \hat V_2^\top \hat V_2 \tilde A_{11}= \tilde A_{11}$ and 
similarly $\hat N_{i, 22} = \tilde N_{i, 11}$ using that $\hat V_2^\top \hat V_2 = I$. Taking \eqref{rightupperzero} into account, the left upper block of \eqref{parttransformedsys} is \begin{align}\nonumber
   & \hat A_{11}\hat P_{11} +  \hat P_{11}\hat A_{11}^\top  + \sum_{i=1}^q   \hat N_{i, 11}\hat P_{11} \hat N_{i, 11}^\top = 0   \\  & \Leftrightarrow   
   \hat P_{11}^{-1}\hat A_{11}^\top + \hat A_{11}\hat P_{11}^{-1} +  \sum_{i=1}^q   \hat P_{11}^{-1}\hat N_{i, 11}\hat P_{11} \hat N_{i, 11}^\top \hat P_{11}^{-1}   = 0.         \label{equationa}      
                                                                                                                                                                                \end{align}
The evaluation of the right upper block yields \begin{align}\nonumber
   & \hat A_{11}\hat P_{12} +  \hat P_{11}\hat A_{21}^\top +\hat P_{12} \hat A_{22}^\top  + \sum_{i=1}^q  \smat {\hat N_{i, 11}}&0\srix 
   \smat {\hat P_{11}}& {\hat P_{12}}\\ {\hat P_{12}^\top}& {\hat P_{22}}\srix \smat {\hat N_{i, 21}^\top}\\  {\hat N_{i, 22}^\top}\srix= 0   \\  & \Leftrightarrow   
  \hat A_{21}^\top =   -\left(\hat P_{11}^{-1}\hat A_{11}\hat P_{12} +\hat P_{11}^{-1}\hat P_{12} \hat A_{22}^\top  + \sum_{i=1}^q  \smat {\hat P_{11}^{-1}\hat N_{i, 11}}&0\srix 
   \smat {\hat P_{11}}& {\hat P_{12}}\\ {\hat P_{12}^\top}& {\hat P_{22}}\srix \smat {\hat N_{i, 21}^\top}\\  {\hat N_{i, 22}^\top}\srix\right).   \label{equationb}
                                                                                                                                                                                \end{align}                                                                                                                                                                                
Finally, the right lower block is given by \begin{align}\label{equationc}
 \hat A_{21}\hat P_{12} + \hat A_{22}\hat P_{22} + \hat P_{12}^\top \hat A_{21}^\top +\hat P_{22} \hat A_{22}^\top  + \sum_{i=1}^q  \smat {\hat N_{i, 21}}&{\hat N_{i, 22}}\srix 
   \smat {\hat P_{11}}& {\hat P_{12}}\\ {\hat P_{12}^\top}& {\hat P_{22}}\srix \smat {\hat N_{i, 21}^\top}\\  {\hat N_{i, 22}^\top}\srix= -\hat R.                                           
                                           \end{align}
We set $\mathbf{\hat P}_{22}= \hat P_{22}- \hat P_{12}^\top \hat P_{11}^{-1} \hat P_{12}$, $\mathbf{\hat N}_{i, 21}=\hat N_{i, 21}- \hat P_{12}^\top \hat P_{11}^{-1} \hat N_{i, 11} $ and insert \eqref{equationb} into \eqref{equationc}
in order to obtain  \begin{align*}
 &\hat A_{22}\mathbf{\hat P}_{22}+\mathbf{\hat P}_{22} \hat A_{22}^\top - \hat P_{12}^\top (\hat A_{11}^\top \hat P_{11}^{-1} +\hat P_{11}^{-1} \hat A_{11})\hat P_{12}   + \sum_{i=1}^q  \smat {\mathbf{\hat N}_{i, 21}}&{\hat N_{i, 22}}\srix 
   \smat {\hat P_{11}}& {\hat P_{12}}\\ {\hat P_{12}^\top}& {\hat P_{22}}\srix \smat {\hat N_{i, 21}^\top}\\  {\hat N_{i, 22}^\top}\srix\\
   & +\sum_{i=1}^q \smat {\hat N_{i, 21}} &  {\hat N_{i, 22}}\srix   
   \smat {\hat P_{11}}& {\hat P_{12}}\\ {\hat P_{12}^\top}& {\hat P_{22}}\srix \smat {-(\hat P_{12}^\top \hat P_{11}^{-1}\hat N_{i, 11})^\top}\\0\srix = -\hat R.                                           
                                           \end{align*}
Using \eqref{equationa} for the above relation leads to \begin{align*}
 &\hat A_{22}\mathbf{\hat P}_{22}+\mathbf{\hat P}_{22} \hat A_{22}^\top   + \sum_{i=1}^q  \smat {\mathbf{\hat N}_{i, 21}}&{\hat N_{i, 22}}\srix 
   \smat {\hat P_{11}}& {\hat P_{12}}\\ {\hat P_{12}^\top}& {\hat P_{22}}\srix \smat {\mathbf{\hat N}_{i, 21}^\top}\\  {\hat N_{i, 22}^\top}\srix = -\hat R.                                           
                                           \end{align*}
 We add and subtract $\sum_{i=1}^q  \hat N_{i, 22} \mathbf{\hat P}_{22} \hat N_{i, 22}^\top$ resulting in \begin{align*}
 &\hat A_{22}\mathbf{\hat P}_{22}+\mathbf{\hat P}_{22} \hat A_{22}^\top   + \sum_{i=1}^q  \hat N_{i, 22} \mathbf{\hat P}_{22} \hat N_{i, 22}^\top+
 \sum_{i=1}^q  \smat {\mathbf{\hat N}_{i, 21}}&{\hat N_{i, 22}}\srix 
   \smat {\hat P_{11}}& {\hat P_{12}}\\ {\hat P_{12}^\top}& {\hat P_{12}^\top \hat P_{11}^{-1} \hat P_{12}}\srix \smat {\mathbf{\hat N}_{i, 21}^\top}\\  {\hat N_{i, 22}^\top}\srix = -\hat R.                                           
                                           \end{align*}
$\smat {\hat P_{11}}& {\hat P_{12}}\\ {\hat P_{12}^\top}& {\hat P_{12}^\top \hat P_{11}^{-1} \hat P_{12}}\srix$ is positive semidefinite since it holds that \begin{align*}
  \smat y^\top & z^\top\srix \smat {\hat P_{11}}& {\hat P_{12}}\\ {\hat P_{12}^\top}& {\hat P_{12}^\top \hat P_{11}^{-1} \hat P_{12}}\srix  \smat y\\z\srix
  =y^\top \hat P_{11} y + 2 y^\top \hat P_{12} z + z^\top \hat P_{12}^\top\hat P_{11}^{-1}\hat P_{12}z = \left\|\hat P_{11}^{\frac{1}{2}} y +  \hat P_{11}^{-\frac{1}{2}}\hat P_{12}z\right\|_2^2\geq 0,
                                                                                                                                                             \end{align*}
where $\smat y\\z\srix$ is an arbitrary vector of suitable dimension. Therefore, we have \begin{align}\nonumber
 \hat A_{22}\mathbf{\hat P}_{22}+\mathbf{\hat P}_{22} \hat A_{22}^\top   + \sum_{i=1}^q  \hat N_{i, 22} \mathbf{\hat P}_{22} \hat N_{i, 22}^\top &=-\left(\hat R + 
 \sum_{i=1}^q  \smat {\mathbf{\hat N}_{i, 21}}&{\hat N_{i, 22}}\srix \smat {\hat P_{11}}& {\hat P_{12}}\\ {\hat P_{12}^\top}& {\hat P_{12}^\top \hat P_{11}^{-1} \hat P_{12}}\srix
 \smat {\mathbf{\hat N}_{i, 21}^\top}\\  {\hat N_{i, 22}^\top}\srix\right)\\&\leq 0            \label{new_red_eq}                            
                                           \end{align}
and $\mathbf{\hat P}_{22}>0$ since it is the inverse of the right lower block of $\smat {\hat P_{11}}& {\hat P_{12}}\\ {\hat P_{12}^\top}& {\hat P_{22}}\srix^{-1}$. By Lemma \ref{lem_stab_plus_zero}, 
this implies $\lambda(I\otimes \hat A_{22} + \hat A_{22}\otimes I+ \sum_{i=1}^{q} \hat N_{i, 22} \otimes \hat N_{i, 22})\subset \overline{\mathbb C_-}$. Let $\hat \Phi_2$ denote the fundamental solution of the 
system with matrices  $(\hat A_{22}, \hat N_{i, 22})$. Moreover, we set $\hat B_2 := \hat V_2^\top B_1$. Then, we can express\begin{align}\label{onesteptilstab}
                 \mathbb E\left\|\hat \Phi(t) B_1\right\|_F^2 = \mathbb E\left\|\hat S^\top(\hat S\hat \Phi(t) \hat S^\top) \hat S B_1\right\|_F^2 = \mathbb E\left\|(\hat S\hat \Phi(t) \hat S^\top) \smat 0\\\hat B_2\srix\right\|_F^2.
                                                                                                                                \end{align}
We partition $\hat S\hat \Phi(t) \hat S^\top=\smat {\hat \Phi_{11}(t)}& {\hat \Phi_{12}(t)}\\ {\hat \Phi_{21}(t)}& {\hat \Phi_{22}(t)}\srix$ and find the associated equation by multiplying the one for $\hat \Phi$ with $\hat S$ from the left 
and $\hat S^\top$ from the right resulting in \begin{align}\label{part22}
\smat {\hat \Phi_{11}}& {\hat \Phi_{12}}\\ {\hat \Phi_{21}}& {\hat \Phi_{22}}\srix=\smat I& 0\\0&I\srix 
+ \int_0^t \smat {\hat A_{11}}& 0\\ {\hat A_{21}}& {\hat A_{22}}\srix \smat {\hat \Phi_{11}}& {\hat \Phi_{12}}\\ {\hat \Phi_{21}}& {\hat \Phi_{22}}\srix ds
+\sum_{i=1}^{q} \int_0^t \smat {\hat N_{i, 11}}& 0\\ {\hat N_{i, 21}}& {\hat N_{i, 22}}\srix \smat {\hat \Phi_{11}}& {\hat \Phi_{12}}\\ {\hat \Phi_{21}}& {\hat \Phi_{22}}\srix dW_i(s).                     
                                              \end{align}
Evaluating the right upper block and subsequently the right lower block of \eqref{part22}, we see that $\hat \Phi_{12}=0$ and $\hat \Phi_{22} =\hat \Phi_{2}$. Therefore, \eqref{onesteptilstab} becomes 
\begin{align}\label{onesteptilstab2}
                 \mathbb E\left\|\hat \Phi(t) B_1\right\|_F^2 = \mathbb E\left\|\hat \Phi_{2}(t)\hat B_2\right\|_F^2.
\end{align}

In addition, we obtain the following rank relation \begin{equation}\begin{aligned}\label{rank_cond}
r_0:&= \rank(\smat B_1 & A_{11}B_1 & \ldots & A_{11}^{r-1}B_1\srix) =  \rank(\smat {\hat S} B_1 & (\hat S A_{11}\hat S^\top){\hat S B_1} & \ldots & (\hat S A_{11}\hat S^\top)^{r-1}{\hat S B}\srix)\\
&= \rank(\smat {\hat B_2} & {\hat A_{22}}{\hat B_2} & \ldots & {\hat A_{22}}^{r-1} {\hat B_2}\srix) = \rank(\smat {\hat B_2} & {\hat A_{22}}{\hat B_2} & \ldots & {\hat A_{22}}^{r_2-1} {\hat B_2}\srix),
                                                    \end{aligned}\end{equation}
where $r_2$ is the number of rows/columns of $\hat A_{22}$.                                                   
If there is no zero eigenvalue of the Kronecker matrix associated to $(\hat A_{22}, \hat N_{i, 22})$, then $\hat \Phi_{2}$ decays exponentially and the claim of this theorem follows by \eqref{onesteptilstab2}. If the projected 
system still has a zero eigenvalue, then by Lemma \ref{lem_zero_eig}, there is $\hat V_{22}\geq 0$, $\hat V_{22}\neq 0$, such that \begin{align*}
       \hat A_{22}^\top{\hat V}_{22}+{\hat V}_{22} \hat A_{22}   + \sum_{i=1}^q  \hat N_{i, 22}^\top {\hat V}_{22} \hat N_{i, 22} = 0.                                                                                                                                                 
                                                                                                                                                       \end{align*}
Now, one can further project down the reduced system with matrices $(\hat A_{22}, \hat B_2, \hat N_{i, 22})$ by the same type of state space transformation as in \eqref{part_trans_red_sys} based on the factor of the 
eigenvalue decomposition of ${\hat V}_{22}$ instead of $\hat S$ and based on \eqref{new_red_eq} instead of \eqref{left_uper_block}. Notice that ${\hat V}_{22}$ cannot have full rank since else we have $\hat B_2 = \hat V_2^\top B_1 = 0$ which, together 
with \eqref{zero_coef2}, implies $B_1 = 0$. One proceeds with this procedure until a mean square asymptotically stable subsystem is achieved. Such a subsystem exists since if one reaches a system of dimension $r_0$, then 
it holds that $r_2 = r_0$ in \eqref{rank_cond}. This local reachability condition combined with \eqref{new_red_eq} is equivalent to mean square asymptotic stability, see \cite[Theorem 3.6.1]{damm}. Since \eqref{onesteptilstab2} is then also obtained with 
the mean square asymptotically stable subsystem, the result follows which completes the proof.

\end{proof}
\begin{remark}
The only structure of the reduced system that was used in the proof of Theorem \ref{thm_red_gram_exisits} is the existence of an equation of the form \eqref{left_uper_block}. 
Therefore, this theorem can be extended to any reduced system for which there exists a matrix $\hat X>0$ such that \begin{align*}
     \hat A \hat X+ \hat X \hat A^\top + \sum_{i=1}^q \hat N_{i} \hat X \hat N_{i}^\top \leq -\hat B \hat B^\top.                                     
                                          \end{align*}
\end{remark}

The following implication of Theorem \ref{thm_red_gram_exisits} shows the square mean asymptotic stability of the ROM \eqref{red_sys_one_sided_BT} is preserved in some extended way.

\begin{kor}\label{kor:MinimalStable} Suppose that the reduced system \eqref{red_stochstate_mul} with matrices $\hat A = A_{11}$, $\hat B = B_1\neq 0$ and $\hat N_{i} =   {N}_{i, 11}$ defined in \eqref{partition}, and
associated to the fundamental solution $\hat \Phi$, is not mean square asymptotically stable. If we further have that $\Lambda_1>0$, then 
there exists $V_0 \in \R^{r \times r_0} $, $r_0<r$, with $V_0^\top V_0 = I$ leading to a  projected system $\hat A_0 = V_0^\top A_{11}V_0$,   $\hat B_0 = V_0^\top B_1 $ and $\hat N_{0, i} =  V_0^\top{N}_{i, 11}V_0$, 
associated to the mean square asymptotically stable fundamental solution $\hat \Phi_0$. Moreover, it holds that 
\[ \hat \Phi(t)B_1 = V_0 \hat{\Phi}_0(t)\hat B_0. \]
\end{kor}
\begin{proof}
As in \eqref{onesteptilstab}, we can write 
$\hat \Phi(t) B_1 = \hat S^\top(\hat S\hat \Phi(t) \hat S^\top) \hat S B_1 $ with the orthogonal matrix $\hat S^\top= \smat {\hat V_1} & {\hat V_2} \srix$. Following the steps of the proof of Theorem
\ref{thm_red_gram_exisits}, we see that $\hat \Phi(t) B_1 = \hat V_2 \hat \Phi_{2}(t) \hat B_2$, where $\hat B_2 = \hat V_2^\top B_1 $ and where $\hat \Phi_{2}$ is the fundamental solution for the system with matrices
$\hat V_2^\top A_{11}\hat V_2$ and $\hat V_2^\top{N}_{i, 11}\hat V_2$. If $\hat \Phi_{2}$ is mean square asymptotically stable, we have that $V_0 = \hat V_2$. Else, by the proof of Theorem
\ref{thm_red_gram_exisits}, the projection procedure can be repeated 
until an mean square asymptotically stable subsystem is achieved. In this case, $V_0$ is the product of matrices like $\hat V_2$.
\end{proof}
Corollary \ref{kor:MinimalStable} shows that the obtained ROM always has a mean square asymptotically stable realization. 
In other words, the procedure described in Section \ref{romsec} produces a ROM that is either mean square asymptotically stable or that can be further reduced to a mean square asymptotically stable system without an 
additional approximation error given that $x_0 = 0$. 
\smallskip

The following corollary will be useful for interpreting error bounds for the approximation error in Section \ref{subsec:ErrBound}.  \smallskip

\begin{kor}\label{gramlesslambda1}
Given the assumptions of Theorem \ref{thm_red_gram_exisits}, we have that \begin{align*}
                                                         \trace(\hat P)\leq \trace(\Lambda_1),
                                                                          \end{align*}
where $\hat P$ is the reachability Gramian of the reduced system with coefficients $\hat A = A_{11}$, $\hat B = B_1$ and $\hat N_{i} =   {N}_{i, 11}$.
\end{kor}
\begin{proof}
As in the proof of Theorem \ref{thm_red_gram_exisits}, three cases need to be considered. Let us first assume that the reduced system is mean square asymptotically stable, i.e., $0\not\in\lambda(\hat K)$. Subtracting 
\eqref{lyap_red_sys} from \eqref{left_uper_block} we see that $\Lambda_1 - \hat P$ satisfies \begin{align}\label{eq_lambda_minus_P}
             A_{11} (\Lambda_1 - \hat P) +  (\Lambda_1 - \hat P) A_{11}^\top + \sum_{i=1}^q N_{i, 11}  (\Lambda_1 - \hat P) N_{i, 11}^\top =  - \sum_{i=1}^q N_{i, 12} \Lambda_2 N_{i, 12}^\top =:-R_2.
                  \end{align}
Now, equation \eqref{eq_lambda_minus_P} is uniquely solvable. According to Section \ref{domsubsec}, this solution is represented by $\mathbb E \int_0^\infty \hat \Phi(s) R_2 \hat \Phi^\top (s) ds\geq 0$. Therefore, we have that 
$\Lambda_1 \geq \hat P$ implying the claim of this corollary. Now, let us study the case of $0\in\lambda(\hat K)$. $B_1=0$ implies that $\hat P = 0$ leading to $\Lambda_1 \geq \hat P$. It remains to consider the case of an 
unstable reduced system with $B_1\neq 0$. We use the arguments of the proof of Theorem \ref{thm_red_gram_exisits} and assume w.l.o.g. that the projected reduced system with matrices
 $(\hat A_{22}, \hat B_2, \hat N_{i, 22})$ and fundamental solution $\hat \Phi_2$ is already mean square asymptotically stable. 
 Else we could project down the reduced system further and the same arguments apply as the ones we use below. Integrating both sides of \eqref{onesteptilstab2} over $[0, \infty)$ and 
 using the definition of the Frobenius norm, we obtain\begin{align*}
                 \trace\bigg(\underbrace{\mathbb E \int_0^\infty\hat \Phi(t) B_1 B_1^\top \hat\Phi^\top(t) dt}_{=\hat P}\bigg)  = 
                 \trace\bigg(\underbrace{\mathbb E\int_0^\infty\hat \Phi_{2}(t)\hat B_2 \hat B_2^\top \hat \Phi_{2}^\top(t) dt}_{=:\hat P_2}\bigg).
\end{align*}
Due to the mean square asymptotic stability, we know that $\hat P_2$ is the unique solution to \begin{align*}
 \hat A_{22}{\hat P}_{2}+{\hat P}_{2} \hat A_{22}^\top   + \sum_{i=1}^q  \hat N_{i, 22} {\hat P}_{2} \hat N_{i, 22}^\top =-  \hat B_2 \hat B_2^\top.                            
                                           \end{align*}
Comparing this equation with \eqref{new_red_eq}, we find that ${\hat P}_{2}\leq \mathbf{\hat P}_{22}= \hat P_{22}- \hat P_{12}^\top \hat P_{11}^{-1} \hat P_{12}\leq \hat P_{22}$. 
We exploit \eqref{part_trans_red_sys} leading to $\trace(\hat P)\leq \trace(\hat P_{22})\leq \trace(\hat P_{22})+\trace(\hat P_{11}) = \trace(\hat S \Lambda_{1} \hat S^\top) = \trace(\Lambda_{1})$.
\end{proof}

\subsection{Error bounds}\label{subsec:ErrBound}
In this subsection, we derive error bounds for the model reduction procedure proposed in Section \ref{romsec}. We begin with an error bound that is general in the sense that it only requires the existence of the Gramians $P$ and $\hat P$ and does not exploit any further structure of the reduced system. Once this general bound is established, 
an error estimate for the choice in \eqref{red_sys_one_sided_BT} is given allowing to identify the scenarios in which this ROM leads to a good approximation. 
The next result characterizes the error in a full state approximation. Notice that we use similar techniques as in \cite{redmannbenner, mliopt}, where output errors were considered. However, we state the following proposition under milder assumptions.
\begin{prop}\label{prop_first_bound}
Suppose that $\Phi$ denotes the fundamental solutions of \eqref{stochstate}, and $\hat \Phi$ denotes the fundamental solutions of \eqref{red_stochstate_mul} obtained by Galerkin projection
using $V \in \R^{n \times r}$ with $V^\top V = I$. Moreover, let $x$ and $\hat x$ represent the solutions to both systems.
If there is a constant $c>0$ such that 
$\mathbb E\left\|\Phi(t) B\right\|^2 $, $\mathbb E\left\|\hat \Phi(t) \hat B\right\|^2 \lesssim \expn^{-c t}$ and if $x_0=0$ and $\hat x_0=0$, 
we have 
\begin{align*}
 \sup_{t\in [0, T]}\mathbb E \left\|x(t) - V\hat x(t) \right\|_2 \leq \left(\trace(P) +   \trace(\hat P) - 2 \trace(P_2 V^\top)\right)^{\frac{1}{2}} \left\| u\right\|_{L^2_T},
\end{align*}
where the matrices $P:= \mathbb E \int_0^\infty \Phi(s)BB^\top \Phi^\top (s) ds$, $\hat P:= \mathbb E \int_0^\infty \hat \Phi(s)\hat B\hat B^\top \hat \Phi^\top (s) ds$, 
$P_2:= \mathbb E \int_0^\infty \Phi(s)B\hat B^\top \hat \Phi^\top (s) ds$ satisfy
 \begin{subequations}\label{all_eq_in_one}
\begin{align}\label{full_gram}
A P + P A^\top+ \sum_{i=1}^{q} N_i P N_{i}^\top &= -B B^\top,\\ \label{red_gram}
\hat A\hat P + \hat P {\hat A}^\top+\sum_{i=1}^{q} \hat N_{i} \hat P {\hat N}_{i}^\top  &= -\hat B \hat B^\top,\\ \label{mixed_gram}
A P_2 + P_2 {\hat A}^\top+\sum_{i=1}^{q} N_{i} P_2 {\hat N}_{i}^\top&= -B\hat B^\top.
\end{align}
\end{subequations}
\end{prop}
\begin{proof}
It can be shown that the solution of \eqref{stochstate} is given by \begin{align*}
x(t) = \Phi(t) x_0 + \int_0^t \Phi(t, s) B u(s) ds,
\end{align*}
see, e.g., \cite{damm}. Setting the initial states in \eqref{stochstate} and \eqref{red_stochstate_mul} equal to zero and using the solution representations for both systems, we obtain by the 
triangle inequality that
\begin{align*}
 \mathbb E \left\|x(t) - V\hat x(t) \right\|_2 &\leq \mathbb E \int_0^t \left\|\left(\Phi(t, s) B - V\hat \Phi(t, s) \hat B\right)u(s) \right\|_2 ds\\  
 &\leq \mathbb E \int_0^t \left\|\Phi(t, s) B - V \hat \Phi(t, s) \hat B\right\|_F \left\|u(s) \right\|_2 ds.
\end{align*}
We apply the inequality of Cauchy-Schwarz and obtain \begin{align*}
 \mathbb E \left\|x(t) - V\hat x(t) \right\|_2 &\leq \left(\mathbb E \int_0^t \left\|\Phi(t, s) B - V \hat \Phi(t, s) \hat B\right\|_F^2 ds\right)^{\frac{1}{2}} \left\| u\right\|_{L^2_t}.
 \end{align*}
The definition of the Frobenius norm and properties of the trace operator yield 
 \begin{align*}
\mathbb E \left\|\Phi(t, s) B - V \hat \Phi(t, s) \hat B\right\|_F^2 = &\trace\left(\mathbb E \left[\Phi(t, s) BB^\top \Phi^\top(t, s)\right]\right)\\ 
&+ \trace\left(V\mathbb E \left[\hat \Phi(t, s) \hat B\hat B^\top \hat \Phi^\top(t, s)V^\top\right]\right)\\
 &-2 \trace\left(\mathbb E \left[\Phi(t, s) B\hat B^\top \hat \Phi^\top(t, s)\right]V^\top\right).
 \end{align*}
Using Corollary \ref{kor_semi_group}, $\Phi(t, s)$ and $\hat \Phi(t, s)$ can be replaced by $\Phi(t-s)$ and $\hat \Phi(t-s)$ above. Writing the resulting trace expressions by the Frobenius norm again, we obtain 
\begin{align*}
 \mathbb E \left\|x(t) - V\hat x(t) \right\|_2 &\leq \left(\mathbb E \int_0^t \left\|\Phi(t - s) B - V \hat \Phi(t - s) \hat B\right\|_F^2 ds\right)^{\frac{1}{2}} \left\| u\right\|_{L^2_t}\\
 &= \left(\mathbb E \int_0^t \left\|\Phi(s) B - V\hat \Phi(s) \hat B\right\|_F^2 ds\right)^{\frac{1}{2}} \left\| u\right\|_{L^2_t}\\
 &\leq \left(\mathbb E \int_0^\infty \left\|\Phi(s) B - V\hat \Phi(s) \hat B\right\|_F^2 ds\right)^{\frac{1}{2}} \left\| u\right\|_{L^2_t}.
\end{align*}
The infinite integral above exists due to the exponential decay of $\Phi B$ and $\hat \Phi \hat B$. Taking the supremum over $[0, T]$, inserting the definition of the Frobenius norm
and exploiting that $V^\top V = I$, we obtain 
\begin{align*}
 \sup_{t\in [0, T]}\mathbb E \left\|x(t) - V\hat x(t) \right\|_2 \leq \left(\trace(P) +   \trace(\hat P) - 2 \trace(P_2 V^\top)\right)^{\frac{1}{2}} \left\| u\right\|_{L^2_T}.
\end{align*}
The infinite integrals $P$, $\hat P$ and $P_2$ satisfy \eqref{all_eq_in_one} due to Lemma \ref{lemdgl} using the exponential decay of $\Phi B$ and $\hat \Phi \hat B$.
\end{proof}
\begin{remark}
Under the assumptions of Proposition \ref{prop_first_bound}, the solutions of \eqref{all_eq_in_one} are not necessarily unique as Example \ref{example2} shows. Uniqueness can be ensured if we further have 
that $\Phi$ and $\hat\Phi$ decay exponentially in the mean square sense.
\end{remark}
Based on the result in Proposition \ref{prop_first_bound}, we now find an error bound for the reduced system introduced in Section \ref{rom_computattion}. Output error bounds for balanced truncation in the same norm based
on different choices of Gramians are proved in \cite{redmannbenner, BTtyp2EB}. The error analysis for the scheme in Section \ref{rom_computattion} is more challenging since less structure than in the case of balanced truncation 
can be exploited which is a method where the reachability and observability Gramian are both diagonal and equal (after a balancing transformation). Moreover, in contrast to balanced truncation, we need 
to discuss the case in which mean square asymptotic stability is not preserved.
\begin{thm}\label{special_error_bound}
Let $x$ be the solution to the mean square asymptotically stable system  \eqref{stochstate} and $\hat x$ the solution to \eqref{red_stochstate_mul} with zero initial states and with 
$\hat A = A_{11}$, $\hat B = B_1$, $\hat N_{i} =   {N}_{i, 11}$ being submatrices of the balanced partition in \eqref{partition}. Let $\Lambda = \diag(\Lambda_1, \Lambda_2)$ be the matrix of ordered 
eigenvalues of the reachability Gramian $P$ with $\Lambda_1 = \diag(\lambda_1, \ldots, \lambda_r)>0$. Let $S^\top =  \smat V & S_{2}\srix$ denote the factor of the associated eigenvalue decomposition of $P$. 
Then, it holds that \begin{align}\label{interprete_bound}
 \sup_{t\in [0, T]}\mathbb E \left\|x(t) - V\hat x(t) \right\|_2 \leq  \left(\trace(\hat P-\Lambda_1) + \trace(\Lambda_2 \mathcal W_0)\right)^{\frac{1}{2}}  \left\| u\right\|_{L^2_T},
\end{align}
where $\hat P$ is the reduced reachability Gramian and \begin{align*}
 \mathcal W_0 = I + 2{A}^\top_{12} Y_2 + \sum_{i=1}^{q} {N}^\top_{i, 12} \left(2Y \smat {N}_{i, 12}\\ {N}_{i, 22}\srix \right).                 
                    \end{align*}
The matrix $Y =  \smat{Y}_{1} & Y_2 \srix$ is defined as the unique solution to 
 \begin{align}\label{yequation}
A_{11}^\top Y+ Y A_b +\sum_{i=1}^{q}  {N}^\top_{i, 11} Y N_{i, b} = -(SV)^\top = -  \smat I & 0 \srix.
    \end{align}   
If it moreover holds that $0\not\in \lambda(I\otimes A_{11} + A_{11}\otimes I+ \sum_{i=1}^{q} N_{i, 11} \otimes N_{i, 11})$, then $\hat Q$ can be introduced as the positive semidefinite 
solution to \begin{align}\label{obs_red_gram}
     A_{11}^\top\hat Q+ \hat Q A_{11}+\sum_{i=1}^{q} N_{i, 11}^\top \hat Q N_{i, 11}  = -I.
                         \end{align}
Hence, the error bound becomes
\begin{align*}
 \sup_{t\in [0, T]}\mathbb E \left\|x(t) - V\hat x(t) \right\|_2 \leq  \left(\trace(\Lambda_2 \mathcal W)\right)^{\frac{1}{2}}  \left\| u\right\|_{L^2_T},
\end{align*}
where the weight is \begin{align*}
 \mathcal W =     I + 2{A}^\top_{12} Y_2 + \sum_{i=1}^{q} {N}^\top_{i, 12} \left(2Y \smat {N}_{i, 12}\\ {N}_{i, 22}\srix - \hat Q N_{i, 12}\right).                
                    \end{align*}
        
\end{thm}
\begin{proof}
Since the original model is asymptotically mean square stable and due to Theorem \ref{thm_red_gram_exisits}, the assumptions of Proposition \ref{prop_first_bound} are met such that we have 
\begin{align*}
 \sup_{t\in [0, T]}\mathbb E \left\|x(t) - V\hat x(t) \right\|_2 \leq \left(\trace(P) +   \trace(\hat P) - 2 \trace(P_2 V^\top)\right)^{\frac{1}{2}} \left\| u\right\|_{L^2_T}.
\end{align*}
Notice that $P$ uniquely solves \eqref{full_gram}. Since the ROM is mean square stable by Proposition \ref{mean_square_stable_red_sys} and due to Lemma \ref{mixed_stab} $P_2$ is also the unique solution to 
\eqref{mixed_gram}. However, there can still be infinitely many other solutions to \eqref{red_gram} besides $\hat P$. 
Using the balanced realization in \eqref{partition}, the error bound then becomes 
\begin{align}\label{firstEB}
 \sup_{t\in [0, T]}\mathbb E \left\|x(t) - V\hat x(t) \right\|_2 \leq  \left(\trace(\Lambda) +   \trace(\hat P) - 2 \trace(S^\top X V^\top)\right)^{\frac{1}{2}}  \left\| u\right\|_{L^2_T},
\end{align}
where $\Lambda$ and $X= S P_2$ uniquely solve \begin{align}\label{full_gram2}
A_b \Lambda + \Lambda A_b^\top+ \sum_{i=1}^{q} N_{i, b} \Lambda N_{i, b}^\top &= -B_b B_b^\top,\\  \label{mixed_gram2}
A_b X + X A_{11}^\top+\sum_{i=1}^{q} N_{i, b} X {N}_{i, 11}^\top &= -B_b B_1^\top.
    \end{align}
By Lemma \ref{mixed_stab}, there is a unique solution to \eqref{yequation} which we can use to rewrite $ \trace(S^\top X V^\top)=  \trace(Y B_b B_1^\top)$. 
Based on the partition \eqref{partition}, we evaluate the first $r$ columns of \eqref{full_gram2} and obtain \begin{align*}
- B_b B_1^\top &=    A_b \smat{\Lambda}_{1}\\ 
0 \srix + \Lambda \smat{A}^\top_{11}\\{A}^\top_{12}\srix+ \sum_{i=1}^{q} N_{i, b} \Lambda \smat{N}^\top_{i, 11}\\{N}^\top_{i, 12}\srix \\    
&=  \smat{A}_{11}\\ {A}_{21}\srix \Lambda_{1} + \smat \Lambda_1 {A}^\top_{11}\\ \Lambda_2 {A}^\top_{12}\srix+ \sum_{i=1}^{q} \left(\smat{N}_{i, 11}\\ 
{N}_{i, 21}\srix  \Lambda_1 {N}^\top_{i, 11} +\smat {N}_{i, 12}\\ 
{N}_{i, 22}\srix \Lambda_2 {N}^\top_{i, 12}\right).     
                                                                                                                            \end{align*}
Inserting this into $\trace(Y B_b B_1^\top)$ yields \begin{align*}
 -\trace(S^\top X V^\top) &=  \trace\left(Y\left[ \smat{A}_{11}\\ {A}_{21}\srix \Lambda_{1} + \smat \Lambda_1 {A}^\top_{11}\\ \Lambda_2 {A}^\top_{12}\srix+ \sum_{i=1}^{q} \left(\smat{N}_{i, 11}\\ 
{N}_{i, 21}\srix  \Lambda_1 {N}^\top_{i, 11} +\smat {N}_{i, 12}\\ 
{N}_{i, 22}\srix \Lambda_2 {N}^\top_{i, 12}\right)\right]\right)  \\
&=   \trace\left(\Lambda_1\left[Y \smat{A}_{11}\\ {A}_{21}\srix+{A}^\top_{11} Y_1 + \sum_{i=1}^{q}  {N}^\top_{i, 11} Y \smat{N}_{i, 11}\\ 
{N}_{i, 21}\srix \right]\right) \\
&\quad + \trace\left(\Lambda_2\left[{A}^\top_{12} Y_2 + \sum_{i=1}^{q} {N}^\top_{i, 12} Y \smat {N}_{i, 12}\\ 
{N}_{i, 22}\srix \right]\right).
                      \end{align*}
The first $r$ columns of \eqref{yequation} yield \begin{align*}
 -\trace(S^\top X V^\top) = -\trace(\Lambda_1)+ \trace\left(\Lambda_2\left[{A}^\top_{12} Y_2 + \sum_{i=1}^{q} {N}^\top_{i, 12} Y \smat {N}_{i, 12}\\ 
{N}_{i, 22}\srix  \right]\right).
\end{align*}
Inserting this into the bound in \eqref{firstEB} leads to \begin{align*}
              \trace(\Lambda) +   \trace(\hat P) - 2 \trace(S^\top X V^\top) = \trace(\hat P-\Lambda_1) + \trace\left(\Lambda_2\left[I + 2{A}^\top_{12} Y_2 + 2\sum_{i=1}^{q} {N}^\top_{i, 12} Y \smat {N}_{i, 12}\\ 
{N}_{i, 22}\srix \right]\right),
                                                          \end{align*}
which proves \eqref{interprete_bound}.

Now let us consider the case where  $0\not\in \lambda(I\otimes A_{11} + A_{11}\otimes I+ \sum_{i=1}^{q} N_{i, 11} \otimes N_{i, 11})$, i.e., the reduced system is mean square asymptotically stable by Proposition \ref{mean_square_stable_red_sys} 
and Lemma \ref{lem_zero_eig}. Therefore, \eqref{obs_red_gram} has a unique positive semidefinite solution $\hat Q$. 
Subtracting the left upper $r\times r$ block of \eqref{full_gram2} from \eqref{red_gram}, we find 
\begin{align*}
 A_{11} (\hat P-\Lambda_1) + (\hat P-\Lambda_1) A_{11}^\top+ \sum_{i=1}^{q} N_{i, 11} (\hat P-\Lambda_1) N_{i, 11}^\top  = \sum_{i=1}^{q} N_{i, 12} \Lambda_2 N_{i, 12}^\top.
\end{align*}
Hence, we have \begin{align*}
\trace(\hat P-\Lambda_1) &= - \trace\left(\left[A_{11}^\top\hat Q+ \hat Q A_{11}+\sum_{i=1}^{q} N_{i, 11}^\top \hat Q N_{i, 11} \right] (\hat P-\Lambda_1)\right)\\
&= - \trace\left(\hat Q \left[A_{11} (\hat P-\Lambda_1) + (\hat P-\Lambda_1) A_{11}^\top+ \sum_{i=1}^{q} N_{i, 11} (\hat P-\Lambda_1) N_{i, 11}^\top \right]\right)\\
&= - \trace\left(\Lambda_2\sum_{i=1}^{q} N_{i, 12}^\top \hat Q N_{i, 12} \right),
               \end{align*}
which concludes the proof of this theorem.
\end{proof} 
Theorem \ref{special_error_bound} is a vital since it shows the relation between the truncated eigenvalues contained in $\Lambda_2$ and the error of the model reduction procedure. By Corollary \ref{gramlesslambda1}, we know that 
$\trace(\hat P-\Lambda_1)\leq 0$ and therefore \eqref{interprete_bound} shows that the error between $x$ and $V\hat x$ is small if $\Lambda_2$ has small diagonal entries. Consequently, the reduced system is accurate 
if only the small eigenvalues of $P$ are neglected. Moreover, this tells us that the reduced order dimension $r$ can be chosen based on the eigenvalues of $P$ since their order is a good indicator for the error. 
Certainly, the error bound representation in Proposition \ref{prop_first_bound} is more suitable for practical computations than the one in Theorem \ref{special_error_bound}. This is because one only needs to solve for 
$\hat P$ and $P_2$ satisfying \eqref{red_gram} and \eqref{mixed_gram} in addition to the Gramian $P$ which is already computed within the model reduction procedure.

\section{Full state approximation for bilinear systems}\label{sec:MORbil}
In this section, we discuss the extension of the proposed results for the class of bilinear systems. We consider the Galerkin projection based model reduction scheme that was studied in Section \ref{rom_computattion} for deterministic bilinear dynamical systems governed by 
\begin{align}\label{bilstate}
 \dot z(t) = Az(t) +  B u(t) + \sum_{i=1}^{m}  N_i z(t) u_i(t),\quad t\geq 0.
\end{align}
Roughly speaking, \eqref{bilstate} is obtained by replacing the white noise processes $\frac{d W_i}{dt}$ in \eqref{stochstate} ($q=m$) by the $i$th component $u_i$ of the control vector $u\in L^2_T$, which 
we assume henceforth  to be deterministic. Transferring the results from the  linear stochastic to the deterministic bilinear case is not trivial, since from the theoretical point of view \eqref{bilstate} and \eqref{stochstate} are very different, since white noise is not a function. However, 
due to the recently shown relation between stochastic and bilinear systems in \cite{h2_bil}, we are able to establish the results of the previous sections for \eqref{bilstate} in a similar manner. Let us assume that the matrix $A$ is Hurwitz, 
i.e., $\lambda(A)\subset \mathbb C_-$. Writing the solution $z= z(\cdot, z_0, B)$  to \eqref{bilstate} dependent on the initial state $z_0$ and the input matrix $B$, $\lambda(A)\subset \mathbb C_-$ implies 
 \begin{align*}
\left\|z(t, z_0, 0)\right\|_2  \lesssim \expn^{-c t}, \quad c>0,                                                                                  
                                                                                    \end{align*}
for all $z_0\in\mathbb R^n$ if $\int_0^\infty \left\| u(s)\right\|_2^2 ds < \infty$, i.e., the homogeneous equation is asymptotically stable with exponential decay, see \cite{h2_bil}. 
If $N_i$ for all $i=1, \ldots, m$ is sufficiently small, $A$ being Hurwitz implies mean square asymptotic stability in the sense of \eqref{assumptionstab}. This can be, e.g., seen by the 
sufficient condition for \eqref{assumptionstab} in \cite[Corollary 3.6.3]{damm}, see also \cite{wonham}. We can now control the matrices $N_i$ by recalling \eqref{bilstate} with $\gamma>0$ resulting in
\begin{align}\label{rescaled_sys}
 \dot z(t)=Az(t)+[\frac{1}{\gamma}B][\gamma u(t)]+\sum_{i=1}^m [\frac{1}{\gamma}N_i] z(t) [\gamma u_i(t)],
            \end{align}
compare also with \cite{bennerdamm,morcondon2005}, where this technique has also been used. 
If $\gamma$ is sufficiently large, \eqref{assumptionstab} can be guaranteed for the pair $(A, \frac{1}{\gamma}N_i)$ which provides the existence of a unique solution to 
\begin{align}\label{reach_gram_bil}
 A P_{\gamma} + P_{\gamma} A^\top + \frac{1}{\gamma^2}\sum_{i=1}^m N_i P_{\gamma} N_i^\top = -\frac{1}{\gamma^2}B B^\top.
\end{align}
According to Section \ref{domsubsec}, $P_{\gamma}$ is the reachability Gramian of the stochastic system \eqref{stochstate} with coefficients $(A, \frac{1}{\gamma}B, \frac{1}{\gamma}N_i)$. 
Choosing $P_\gamma$ for $\gamma = 1$ as a reachability Gramian in the context of model reduction for bilinear systems was first proposed in \cite{typeIBT} and, e.g., further investigated in \cite{bennerdamm}. 
By \cite{h2_bil} we know that $P_\gamma$ 
takes a similar role as in the stochastic case (compare with \eqref{diffreachjanein}), i.e., it characterizes redundant information in \eqref{bilstate} by \begin{align}\label{reach_est_bil}
  \sup_{t\in [0, T]}\left\vert\langle z(t, 0, B), p_{\gamma, k}\rangle_2\right\vert \leq \lambda_{\gamma, k}^{\frac{1}{2}} \exp\left\{0.5\gamma^2\left\|u^{0}\right\|_{L^2_T}^2\right\} \gamma\left\|u\right\|_{L^2_T},          
                    \end{align}
where $(p_{\gamma, k})$ is an orthonormal basis of eigenvector of $P_{\gamma}$ with associated eigenvalues $(\lambda_{\gamma, k})$ and $u^0$ the vector of controls entering the bilinear part of the equation, i.e., 
\begin{align}\label{uzero}
u^{0}=(u_1^{0}\ u_2^{0}\, \ldots\, u_m^{0})^\top\quad \text{with}\quad u_i^{0} \equiv \begin{cases}
  0,  & \text{if }N_i = 0\\
  u_i, & \text{else}.
\end{cases}
\end{align}
Therefore, the eigenspaces of $P_{\gamma}$ corresponding to the zero eigenvalues are irrelevant for the system dynamics. Moreover, assuming that the control energy is sufficiently small, 
\eqref{reach_est_bil} tells us that $z(\cdot, 0, B)$ is small in the direction of $p_{\gamma, k}$ if $\lambda_{\gamma, k}$ is small. Therefore, these eigenspaces can also be seen as less relevant in \eqref{bilstate} and can 
hence be removed leading to ROMs. A somehow different way of characterizing unimportant states in a bilinear equation was discussed in \cite{bennerdamm, graymesko}, where local estimates for 
the reachability energy based on $P_{\gamma}$, $\gamma =1$ have been shown.
\begin{remark}
So far, we observed some essential differences between stochastic and bilinear systems. System \eqref{bilstate} only requires $A$ to be Hurwitz instead of \eqref{assumptionstab}. On the other hand, we consider a 
family of Gramians for the bilinear case depending on $\gamma$ rather than a fixed Gramian.
Although the characterization of irrelevant states are similar in both cases, the exponential in \eqref{reach_est_bil} indicates that we need a certain 
smallness assumption on $u^0$ and $\gamma$ in order to make our arguments valid.
\end{remark}
The above considerations motivate to conduct the same reduced order modeling procedure as explained in Section \ref{rom_computattion}. We introduce the eigenvalue decomposition of 
\begin{align*}
          P_\gamma= S_\gamma^\top \smat{\Lambda}_{\gamma, 1}& 0\\ 0 &{\Lambda}_{\gamma, 2}\srix  S_\gamma,                                    
                                             \end{align*}
where $\Lambda_{\gamma, 1}> 0$ contains the large and $\Lambda_{\gamma, 2}$ the small ordered eigenvalues of $P_\gamma$. Using the partition
\begin{align}\label{partition_bil}
S_\gamma{A}S_\gamma^\top= \smat{A}_{11}^{(\gamma)}&{A}_{12}^{(\gamma)}\\ 
{A}_{21}^{(\gamma)}&{A}_{22}^{(\gamma)}\srix,\quad S_\gamma {B} = \smat{B}_1^{(\gamma)}\\ {B}_2^{(\gamma)}\srix,\quad  S_\gamma{N_i}S_\gamma^\top= \smat{N}_{i, 11}^{(\gamma)}&{N}_{i, 12}^{(\gamma)}\\ 
{N}_{i, 21}^{(\gamma)}&{N}_{i, 22}^{(\gamma)}\srix\end{align}
the eigenvectors associated to small eigenvalues of $P_\gamma$ are then truncated, resulting in the reduced model \begin{align}\label{bilstate_red}
 \dot {\hat z}_\gamma(t) = {A}_{11}^{(\gamma)}\hat z_\gamma(t) +  {B}_1^{(\gamma)} u(t) + \sum_{i=1}^{m} {N}_{i, 11}^{(\gamma)} \hat z_\gamma(t) u_i(t),\quad t\geq 0.
\end{align}
The properties of \eqref{bilstate_red} can now be immediately transferred from the considerations in the stochastic case. By Proposition \ref{mean_square_stable_red_sys}, we have \begin{align}\label{red_kron_eig_bil}
\lambda(I\otimes A_{11}^{(\gamma)} + A_{11}^{(\gamma)}\otimes I+ \frac{1}{\gamma^2}\sum_{i=1}^{m} N_{i, 11}^{(\gamma)} \otimes N_{i, 11}^{(\gamma)})\subset \overline{\mathbb C_-}.                                                                                                                                                                          
                                                                                                                                                                        \end{align}
Example \ref{example1} shows that eigenvalues on the imaginary axis can occur in \eqref{red_kron_eig_bil}, but they can be excluded by Lemma \ref{lem_zero_eig} if       
$0\not\in \lambda(I\otimes A_{11}^{(\gamma)} + A_{11}^{(\gamma)}\otimes I+ \frac{1}{\gamma^2} \sum_{i=1}^{m} N_{i, 11}^{(\gamma)} \otimes N_{i, 11}^{(\gamma)})$. However, even though there is a zero eigenvalue, the existence of the reduced 
order Gramian can be guaranteed using the arguments of Theorem \ref{thm_red_gram_exisits}.\smallskip

In order to keep the discussion around the error bound for bilinear systems short, we do not discuss the scenario of a zero eigenvalue \eqref{red_kron_eig_bil} in detail. Therefore, let us exclude this case below. 
 We can now transfer the result of Proposition \ref{prop_first_bound} to the bilinear case  by the results of \cite{h2_bil}.
  \begin{prop}\label{prop_first_bound_bil}
Let $z$ be the solution to \eqref{bilstate} with $\lambda(A)\subset \mathbb C_-$ and let $\hat z_\gamma$ represent the solution to \eqref{bilstate_red}. Moreover, let 
$\gamma>0$ such that \begin{align*}
\lambda(I\otimes A + A\otimes I+ \frac{1}{\gamma^2}\sum_{i=1}^{m} N_{i} \otimes N_{i})\subset {\mathbb C_-}                                                                                                                                                                       
                                                                                                                                                                        \end{align*}
and that the ROM coefficients satisfy   \begin{align*}
           0\not\in \lambda(I\otimes A_{11}^{(\gamma)} + A_{11}^{(\gamma)}\otimes I+ \frac{1}{\gamma^2} \sum_{i=1}^{m} N_{i, 11}^{(\gamma)} \otimes N_{i, 11}^{(\gamma)}).                                              
                                                        \end{align*}                                                                                                            
Given zero initial states to both equations and $V_\gamma \in \R^{n \times r}$ being the first $r$ columns of the factor $S^\top_\gamma$ of the eigenvalue decomposition of $P_\gamma$
(unique solution to \eqref{reach_gram_bil}), we have \begin{align*}
 \sup_{t\in [0, T]}\left\|z(t) - V_\gamma\hat z_\gamma(t) \right\|_2 \leq \left(\trace(P_\gamma) +   \trace(\hat P_\gamma) - 2 \trace(P_{_\gamma, 2} V_\gamma^\top)\right)^{\frac{1}{2}} 
 \exp\left\{0.5\gamma^2\left\|u^{0}\right\|_{L^2_T}^2\right\} \gamma\left\|u\right\|_{L^2_T},
\end{align*}
where $P_{\gamma, 2}$ and $\hat P_{\gamma}$ are the unique solutions to
\begin{align*}
A P_{\gamma, 2} + P_{\gamma, 2} {{A}_{11}^{(\gamma)}}^\top+\frac{1}{\gamma^2} \sum_{i=1}^{m} N_{i} P_{\gamma, 2} {{N}_{i, 11}^{(\gamma)}}^\top &= -\frac{1}{\gamma^2} B{{B}_1^{(\gamma)}}^\top,\\
{A}_{11}^{(\gamma)} \hat P_{\gamma} + \hat P_{\gamma} {{A}_{11}^{(\gamma)}}^\top+\frac{1}{\gamma^2} \sum_{i=1}^{m} {N}_{i, 11}^{(\gamma)} \hat P_{\gamma} {{N}_{i, 11}^{(\gamma)}}^\top &= 
-\frac{1}{\gamma^2} {B}_1^{(\gamma)} {{B}_1^{(\gamma)}}^\top.\end{align*}
\end{prop}
 \begin{proof}
Given the assumptions $P_\gamma$, $P_{\gamma, 2}$ and $\hat P_{\gamma}$ exist. The result is then a direct consequence of Corollary 4.3 in \cite{h2_bil}.
\end{proof}
\begin{thm}\label{error_bound_bil_case}
Under the assumptions of Proposition \ref{prop_first_bound_bil}, we have 
\begin{align}\label{bound_bil}
 \sup_{t\in [0, T]} \left\|z(t) - V_\gamma \hat z_\gamma(t) \right\|_2 \leq  \left(\trace(\Lambda_{\gamma, 2}\mathcal W_\gamma)\right)^{\frac{1}{2}} 
 \exp\left\{0.5\gamma^2\left\|u^{0}\right\|_{L^2_T}^2\right\} \gamma\left\|u\right\|_{L^2_T},
\end{align}
where the weight is \begin{align*}
 \mathcal W_\gamma =     I + 2 {A_{12}^{(\gamma)}}^\top Y_{\gamma, 2} 
 + \frac{1}{\gamma^2} \sum_{i=1}^{m} {N_{i, 12}^{(\gamma)}}^\top \left(2Y_\gamma \smat {N}_{i, 12}^{(\gamma)}\\ {N}_{i, 22}^{(\gamma)}\srix - \hat Q_\gamma N_{i, 12}^{(\gamma)}\right).                
                    \end{align*}
Above, $Y_\gamma =  \smat{Y}_{\gamma, 1} & Y_{\gamma, 2} \srix$ and $\hat Q_\gamma$ are defined as the unique solutions to 
 \begin{align*}
{{A}_{11}^{(\gamma)}}^\top Y_\gamma+ Y_\gamma A_b^{(\gamma)} +\frac{1}{\gamma^2} \sum_{i=1}^{m}  {{N}_{i, 11}^{(\gamma)}}^\top Y_\gamma N_{i, b}^{(\gamma)} & = -  \smat I & 0 \srix,\\
    {{A}_{11}^{(\gamma)}}^\top\hat Q_\gamma+ \hat Q_\gamma {A}_{11}^{(\gamma)}+\frac{1}{\gamma^2} \sum_{i=1}^{m} {{N}_{i, 11}^{(\gamma)}}^\top \hat Q_\gamma {N}_{i, 11}^{(\gamma)}  &= -I,
                         \end{align*}    
where we set $A_b^{(\gamma)}:= S_\gamma{A}S_\gamma^\top$ and $N_{i, b}^{(\gamma)}:= S_\gamma{N_i}S_\gamma^\top$.
\end{thm}
\begin{proof}
 The result directly follows from the proof of Theorem \ref{special_error_bound} in which $B$ and $N_i$ need to be replaced by $\frac{1}{\gamma}B$ and $\frac{1}{\gamma}N_i$.
\end{proof}
As in the stochastic framework, we can conclude that truncating the small eigenvalues of $P_\gamma$ leads to small diagonal entries of $\Lambda_{\gamma, 2}$ and hence to a small error in the dimension reduction 
according to Theorem \ref{error_bound_bil_case} given that the exponential in \eqref{bound_bil} is not too dominant. Therefore, the eigenvalues of $P_\gamma$ can be used as a criterion to determine a 
suitable reduced order dimension $r$.

\section{Numerical experiments}\label{sec:NumExp}

In this section, we test the efficiency of the proposed method (see Sections \ref{rom_computattion} and \ref{sec:MORbil}), denoted here by \texttt{OS}, in some numerical examples.  
We compare the results with the ones obtained by applying the standard balanced truncation method for a full state approximation, denoted here by \texttt{BT} (see, e.g., \cite{bennerdamm} for the bilinear and 
\cite{redmannbenner} for the stochastic case). All the simulations are done on a CPU 2.6 GHz \intel~\coreifive, 8 GB 1600 MHz DDR3, \matlab~9.1.0.441655 (R2016b).

For this study, we consider a standard test example representing a  $2$D boundary controlled heat transfer system; see, e.g., \cite{bennerdamm}. 
Its dynamics is governed by the heat equation subject to Dirichlet and Robin boundary conditions, i.e., the following boundary value problem
\begin{align*}
\begin{array}{rll}
\partial_{t}\,x &= \Delta x, & \text{in} \, (0,1)\times(0,1), \\
n\cdot \nabla x &= 0.8u_1x & \text{on} \, \Gamma_{1},\\
x &= u_2  , & \text{on} \, \Gamma_{2},\\
x &= 0, & \text{on} \, \Gamma_{3}, \Gamma_{4}, \\
\end{array}
\end{align*}
where $\Gamma_1 = \{0\}\times (0,1)$, $\Gamma_2 = (0,1)\times\{0\}$, $\Gamma_3 = \{1\}\times (0,1)$ and $\Gamma_4 = (0,1)\times\{1\}$. 
In this system, there are two source terms, namely $u_1$ and $u_2$, which are applied at the boundaries $\Gamma_1$ and $\Gamma_2$, respectively. 
A semi-discretization in space using finite differences with $k = 20$ grid points results in a control system of dimension $n = 400$ of the form
\begin{equation}\label{eq:NumExample}
 \dot{x} = Ax(t) + Nx(t)u_1(t) + Bu_2(t). 
\end{equation}
We refer to \cite{bennerdamm} for more details on the matrices in \eqref{eq:NumExample}.

\subsection{Stochastic example}
First, we consider that the boundary $\Gamma_{1}$ is a perturbed by noise, i.e., $u_1 = \frac{dW}{dt}$ with $W$ being a standard Wiener process. Hence the resulting dynamical stochastic system is of  the form
\[  dx(t) = [Ax(t) +  B u_2(t)]dt +  N x(t) dW(t),\quad t\geq 0. \]

\begin{figure}[tbh]
	\newlength\wex
	\newlength\hex
	\setlength{\wex}{.7\textwidth}
	\setlength{\hex}{0.3\textwidth}
	\begin{center}	
	\input{stochastic_os_vs_bt_SVD_decay_Hsvd_decay.tex}	
	\end{center}
	\caption{Decay of singular values $\sigma_k$: the blue curve corresponds to $\eig{P}$. The red curve corresponds to $\sqrt{\eig{PQ}}$.}\label{fig:SV}
\end{figure}
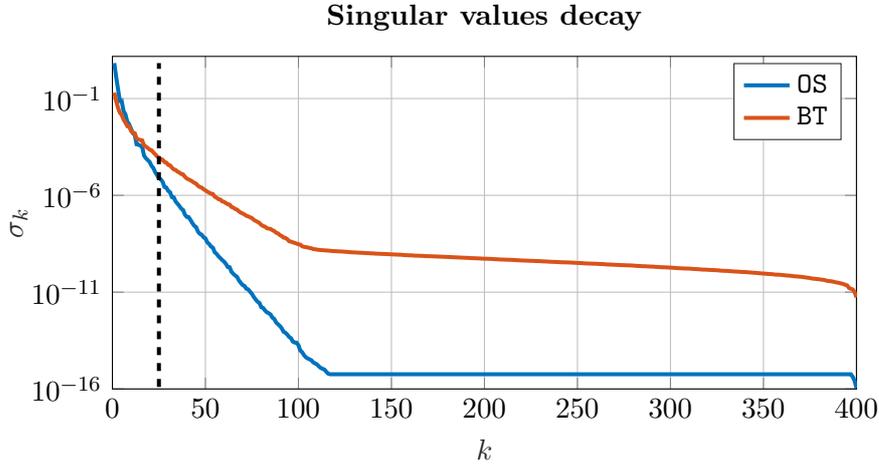

In order to apply \texttt{BT}, we additionally need to compute the observability Gramian $Q$, which satisfies the following Lyapunov equation 
\begin{equation}\label{eq:ObsLyap}
A^\top Q+QA +  \sum_{i=1}^q N_i^\top Q N_i = -I
\end{equation}
with $q=1$ and $N_1 = N$. 
This method was studied in detail in \cite{redmannbenner}. However, solving \eqref{eq:ObsLyap} leads to much higher computational cost especially due to the full-rank right hand side which does not allow the usage 
of low-rank solvers.
Figure \ref{fig:SV} depicts the decay of the eigenvalues/singular values of $P$ as well as the decay of square root of the eigenvalues of $PQ$ (Hankel singular values). As shown in Theorem \ref{special_error_bound}, 
the eigenvalues of $P$ play an important role in the error bound for \texttt{OS} and provide an intuition for the expected error. 
Similarly, the Hankel singular values are also associated the error bound for \texttt{BT}, see \cite{redmannbenner}. The decay of both curves in Figure \ref{fig:SV} 
indicates that a small reduction error can already be achieved for 
for small $r$.

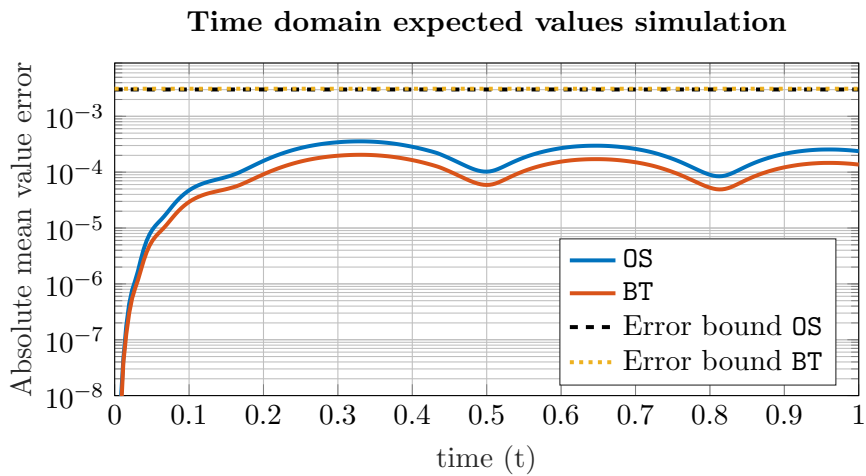
\begin{figure}[tbh]
	\setlength{\wex}{0.7\textwidth}
	\setlength{\hex}{0.3\textwidth}
	\begin{center}
		\input{stochastic_os_vs_bt_r25_nt256_NS1e5.tex}	
		\caption{Stochastic simulation: mean error between the original model and the ROMs for the input $u_2(t) = u(t)= e^{-\frac{1}{2}t}\sin(10t)$.}\label{fig:stoch_abserror}
	\end{center}
\end{figure}

For this example, we compute reduced systems of order $r=25$ for both \texttt{OS} and \texttt{BT}. 
As a next step, we compare the quality of the reduced-order systems by simulating their responses for the input $u_2(t) = u(t) = e^{-\frac{1}{2}t}\sin(10t)$. 
To determine the transient response, we apply a semi-implicit Euler-Maruyama scheme with step size $h = 1/256$ and simulate the original system and the reduced-order models in the time interval $[0,1]$. 
Additionally, those simulations are done using $10^5$ samples. The mean error between the original and the reduced models are depicted in Figure \ref{fig:stoch_abserror} as well as the error bounds 
from Proposition \ref{prop_first_bound}. Table \ref{table:stoch_tb} presents the numerical values for the error bounds and max mean error for both methods. 
We notice that both reduced models are able to follow the behavior of the original system. Furthermore, this figure shows that the two methods, \texttt{BT} and \texttt{OS},  
provide very similar quality reduced models in terms of the magnitude of the error, an observation we also made with other test examples.
However, we note that \texttt{BT} is a numerically more expensive method, since one needs to additionally solve for $Q$. 

\begin{table}[!tb]
	\centering
	\begin{tabular}{|l|l|l|}\hline
		Method & Error bound & max mean error \\ \hline
		\texttt{OS} &      $5.46 \cdot 10^{-3}$      &           $3.56\cdot 10^{-4}$                    \\ \hline
		\texttt{BT}&   $5.63\cdot 10^{-3}$         &           $2.05\cdot 10^{-4}$        \\ \hline            
	\end{tabular}
\caption{Stochastic example: Error bounds and the max value of the mean error  for \texttt{OS} and \texttt{BT} for the simulation presented in Figure \ref{fig:stoch_abserror}.}\label{table:stoch_tb}
\end{table}

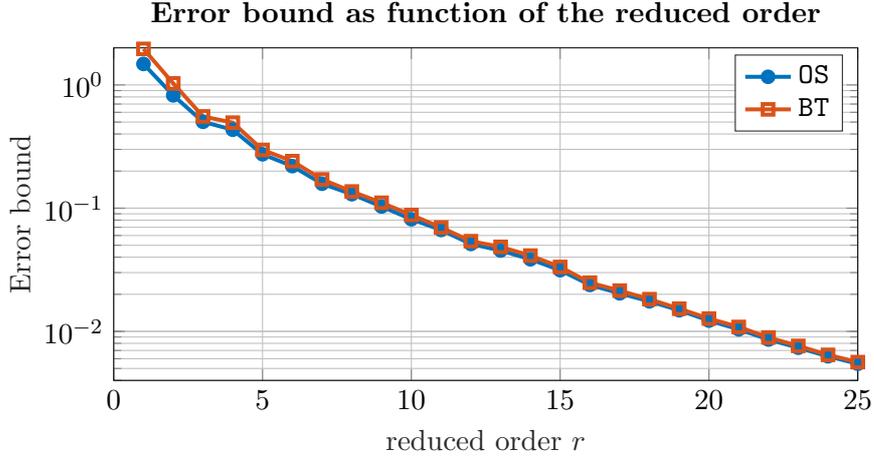
\begin{figure}[tb]
	\setlength{\wex}{0.7\textwidth}
	\setlength{\hex}{0.3\textwidth}
	\begin{center}
		\input{normH2decay_os_vs_bt.tex}	
		\caption{Decay input-independent part of error bound for \texttt{OS} and \texttt{BT} computed for  different orders $r =1, \dots, 25$. }\label{fig:H2error_decay}
	\end{center}
\end{figure}

Additionally, for different reduced orders varying in in the range $r = 1,\dots, 25$,  the input-independent part of the error bound given in Proposition \ref{prop_first_bound} is computed in Figure \ref{fig:H2error_decay}, i.e., 
for each reduced order $r$ we plot the value
\[\cE(r) = \left(\trace(P) +   \trace(\hat P(r)V(r)^\top V(r)) - 2 \trace(P_2(r) V(r)^\top)\right)^{\frac{1}{2}}, \]
where $V(r)$ is the reduced basis of order $r$, and $\hat P(r)$,  $P_2(r)$ are the solutions of \eqref{red_gram} and \eqref{mixed_gram}. 
Notice that we added $V(r)^\top V(r)$ in the second summand of the error bound since $V(r)^\top V(r)\neq I$ for \texttt{BT}. 
As expected, the bound decays if the reduced order is increased for both \texttt{OS} and \texttt{BT}.

\subsection{Bilinear example}
As our second numerical example, we consider the heat transfer system in \eqref{eq:NumExample} with $u_2 = u_1 = u$. 
As a consequence, this leads to a bilinear system having only one input. 
For this example, we need to solve the Lyapunov equation in \eqref{reach_gram_bil}. To this aim, we set $\gamma = 1$ leading to the same reachability and observability Gramians as for the stochastic example. 
Hence, Figure \ref{fig:SV} also gives the decay of singular values for \texttt{OS} and \texttt{BT}. Similarly, Figure \ref{fig:H2error_decay} shows the decay of the input-independent part of the error bound from
Proposition \ref{prop_first_bound_bil}.

\begin{figure}[tb]
	\setlength{\wex}{0.7\textwidth}
	\setlength{\hex}{0.3\textwidth}
	\begin{center}
		\input{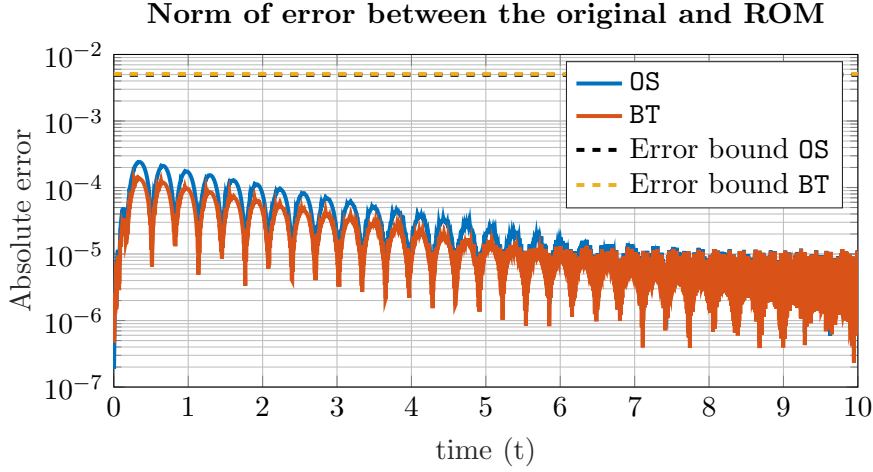}	
		\caption{Bilinear simulation: absolute error between the original model and the ROMs for the input $u(t)= e^{-\frac{1}{2}t}\sin(10t)$.}\label{fig:bil_abserror}
	\end{center}
\end{figure}

\begin{figure}[tb]
	\setlength{\wex}{0.7\textwidth}
	\setlength{\hex}{0.3\textwidth}
	\begin{center}
		\input{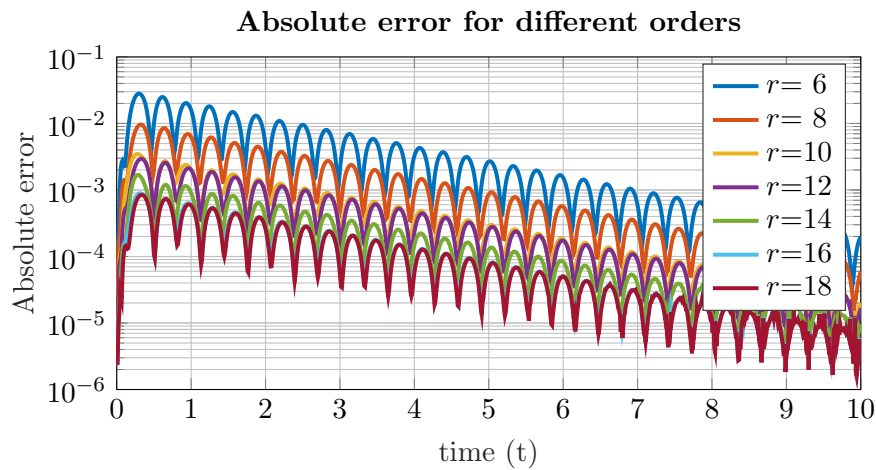}	
		\caption{Bilinear simulation: time domain error for different reduced orders using \texttt{OS}.}\label{fig:bil_error_vs_order}
	\end{center}
\end{figure}

As in the previous example, we obtain reduced systems of order $r=25$ by using \texttt{OS} and \texttt{BT} and compare their quality by
simulating their responses for the input $u(t) = e^{-\frac{1}{2}t}\sin(10t)$. To determine the transient response, we use the \matlab solver \texttt{ode45} to  simulate the original system and the reduced-order models
in the time interval $[0,10]$. The results are depicted in Figure \ref{fig:bil_abserror}. Table \ref{fig:bil_abserror} presents the numerical values for the error bounds and max error for both methods.
Similar to the stochastic example, we notice that the two methods, \texttt{BT} and \texttt{OS},  provide very similar quality reduced models in terms of the magnitude of the error.
Once again, we note that \texttt{BT} is a computationally more expensive method, since one needs the solution to the additional Lyapunov equation in \eqref{eq:ObsLyap}. Finally, Figure \ref{fig:bil_error_vs_order} 
shows the simulation of the error for reduced models obtained by \texttt{OS} with different orders. As expected, the error decays once the order is increased. 
\begin{table}[!tb]
	\centering
	\begin{tabular}{|l|l|l|}\hline
		Method & Error bound & $L_{\infty}$ error \\ \hline
		\texttt{OS} &      $5.46 \cdot 10^{-3}$      &           $3.56\cdot 10^{-4}$                    \\ \hline
		\texttt{BT}&   $5.63\cdot 10^{-3}$         &           $2.05\cdot 10^{-4}$        \\ \hline            
	\end{tabular}\label{table:bil_tb}
	\caption{Bilinear example: Error bounds and the $L_\infty$ error  for \texttt{OS} and \texttt{BT} for the simulation presented in Figure \ref{fig:bil_abserror}. }
\end{table}

\appendix
\section{Matrix differential equations and their solutions}\label{appA}
\begin{lem}\label{lemdgl}
Let $\Phi$ be the fundamental solution of \eqref{stochstate} defined in \eqref{funddef} and let $\hat \Phi$ be the one of system \eqref{red_stochstate_mul}. 
Suppose that $B$ and $\hat B$ are matrices of suitable dimension. Then,  
the $\mathbb R^{n\times r}$-valued function $\mathbb E \left[\Phi(t, s) B {\hat B}^\top \hat{\Phi}^\top(t, s)\right]$, $t\geq s$, satisfies \begin{align}\label{dglmixed}
X(t)= B {\hat B}^\top + \int_s^t  A X(\tau)d\tau + \int_s^t X(\tau) {\hat A}^\top d\tau + \sum_{i = 1}^{q} \int_s^t N_i X(\tau) {\hat N}_i^\top d\tau.
\end{align}
\end{lem}
\begin{proof}
The result is a direct consequence of \cite[Proposition 4.4]{redmannbenner} or \cite[Lemma 2.1]{mliopt}.
\end{proof}
\begin{kor}\label{kor_semi_group}
Given the assumptions in Lemma \ref{lemdgl}, we find that \begin{align}\label{semi_group_prop}
\mathbb E \left[\Phi(t, s) B {\hat B}^\top \hat{\Phi}^\top(t, s)\right] = \mathbb E \left[\Phi(t-s) B {\hat B}^\top \hat{\Phi}^\top(t-s)\right], \quad t\geq s.
\end{align}
\end{kor}
\begin{proof}
Setting $X(t):=\mathbb E \left[\Phi(t) B {\hat B}^\top \hat{\Phi}^\top(t)\right]$, by Lemma \ref{lemdgl} we find that
\begin{align*}
X(t-s)= B {\hat B}^\top + \int_0^{t-s}  A X(\tau)d\tau + \int_0^{t-s} X(\tau) {\hat A}^\top d\tau + \sum_{i = 1}^{q} \int_0^{t-s} N_i X(\tau) {\hat N}_i^\top d\tau.
\end{align*}
Setting $v = \tau + s$, by substitution, we see that \begin{align*}
X(t-s)= B {\hat B}^\top + \int_s^{t}  A X(v-s)dv + \int_s^{t} X(v-s) {\hat A}^\top dv + \sum_{i = 1}^{q} \int_s^{t} N_i X(v-s) {\hat N}_i^\top dv.
\end{align*}
Consequently, both sides of \eqref{semi_group_prop} satisfy \eqref{dglmixed}. Therefore, they are equal.
\end{proof}

\bibliographystyle{plain}

\end{document}

%% file: stochastic_os_vs_bt_SVD_decay_Hsvd_decay.tex
%
%
\definecolor{mycolor1}{rgb}{0.00000,0.44700,0.74100}%
\definecolor{mycolor2}{rgb}{0.85000,0.32500,0.09800}%
\begin{tikzpicture}

\begin{axis}[%
width=0.951\wex,
height=\hex,
at={(0\wex,0\hex)},
scale only axis,
xmin=0,
xmax=400,
xlabel style={font=\color{white!15!black}},
xlabel={$k$},
ymode=log,
ymin=1e-16,
ymax=15,
yminorticks=true,
ylabel style={font=\color{white!15!black}},
ylabel={$\sigma_k$},
axis background/.style={fill=white},
title style={font=\bfseries},
title={Singular values decay},
legend style={legend cell align=left, align=left, draw=white!15!black},
xmajorgrids,
ymajorgrids,
yminorgrids
]
\addplot [color=mycolor1, line width=1.5pt]
  table[row sep=crcr]{%
1	6.76440006564327\\
2	1.3287267974021\\
3	0.314058265203795\\
4	0.0733489836300175\\
5	0.0725772820467477\\
6	0.0224522833169765\\
7	0.0153257547722421\\
8	0.00635015683427786\\
9	0.0046990893828586\\
10	0.00275137768566084\\
11	0.00161979181375633\\
12	0.00136337061601202\\
13	0.000428792601067511\\
14	0.000413065204526199\\
15	0.00036685812684474\\
16	0.000326066195811229\\
17	0.000113066764302002\\
18	8.26814445864444e-05\\
19	6.58617639094522e-05\\
20	5.33259690289805e-05\\
21	3.38206671547712e-05\\
22	2.57384967701659e-05\\
23	1.59008298771336e-05\\
24	1.10705099089441e-05\\
25	7.61061034118793e-06\\
26	6.70598178344214e-06\\
27	5.02188720457399e-06\\
28	2.70636441291298e-06\\
29	2.52353300165658e-06\\
30	1.68725307026553e-06\\
31	1.45568941681086e-06\\
32	8.09594114964378e-07\\
33	7.78887634679756e-07\\
34	4.82041529478094e-07\\
35	3.70894959059838e-07\\
36	2.82550573866366e-07\\
37	2.05808856151608e-07\\
38	1.74874229325545e-07\\
39	1.00000422567225e-07\\
40	7.69588319105811e-08\\
41	7.4813335510989e-08\\
42	5.31275916156835e-08\\
43	3.68462696099188e-08\\
44	2.66550089642099e-08\\
45	2.28535081509574e-08\\
46	1.33618489851761e-08\\
47	1.03225818425708e-08\\
48	9.2745089842486e-09\\
49	6.78466662366432e-09\\
50	6.01988689712205e-09\\
51	3.72904734456468e-09\\
52	3.07975358594923e-09\\
53	2.27217033990669e-09\\
54	1.99193425160427e-09\\
55	1.15008221946512e-09\\
56	9.50356484777702e-10\\
57	7.12789619504168e-10\\
58	5.7576596405233e-10\\
59	4.57096102121496e-10\\
60	3.78109504610548e-10\\
61	2.87864684498378e-10\\
62	1.99352990960932e-10\\
63	1.94136166364251e-10\\
64	1.12296531947841e-10\\
65	8.8612946016992e-11\\
66	7.22469432470188e-11\\
67	6.17787959274262e-11\\
68	4.01583592448095e-11\\
69	3.03235662129786e-11\\
70	2.32401312506643e-11\\
71	2.09512034446309e-11\\
72	1.49750527082352e-11\\
73	1.38866883551214e-11\\
74	1.05173278738387e-11\\
75	8.14877362365295e-12\\
76	5.64050240002851e-12\\
77	4.19312946439268e-12\\
78	3.22454695748081e-12\\
79	2.8072714623874e-12\\
80	1.67286103732119e-12\\
81	1.54101501497374e-12\\
82	1.14013270598606e-12\\
83	9.18030999049552e-13\\
84	8.21490830555491e-13\\
85	7.00718023931968e-13\\
86	4.84323823408348e-13\\
87	4.11365399796849e-13\\
88	2.7181461276204e-13\\
89	2.2000643121379e-13\\
90	1.81362111015596e-13\\
91	1.35737510193231e-13\\
92	9.0042364647813e-14\\
93	7.99499261526852e-14\\
94	6.45671512623112e-14\\
95	4.85285619001213e-14\\
96	4.19266039440221e-14\\
97	2.7908976788479e-14\\
98	2.50531609728912e-14\\
99	2.44901415235127e-14\\
100	1.83421621517831e-14\\
101	9.63314512634315e-15\\
102	7.95451129926849e-15\\
103	5.57282595223008e-15\\
104	4.32571658610306e-15\\
105	3.76831989288562e-15\\
106	3.37181563192055e-15\\
107	3.03666461381695e-15\\
108	2.34436681798685e-15\\
109	2.09902835145827e-15\\
110	1.84671192708019e-15\\
111	1.46997723766665e-15\\
112	1.22052577586463e-15\\
113	1.10172966147024e-15\\
114	9.99726473978139e-16\\
115	7.71518372656475e-16\\
116	6.15934898679024e-16\\
117	5.74737336103367e-16\\
118	5.74737336103367e-16\\
119	5.74737336103367e-16\\
120	5.74737336103367e-16\\
121	5.74737336103367e-16\\
122	5.74737336103367e-16\\
123	5.74737336103367e-16\\
124	5.74737336103367e-16\\
125	5.74737336103367e-16\\
126	5.74737336103367e-16\\
127	5.74737336103367e-16\\
128	5.74737336103367e-16\\
129	5.74737336103367e-16\\
130	5.74737336103367e-16\\
131	5.74737336103367e-16\\
132	5.74737336103367e-16\\
133	5.74737336103367e-16\\
134	5.74737336103367e-16\\
135	5.74737336103367e-16\\
136	5.74737336103367e-16\\
137	5.74737336103367e-16\\
138	5.74737336103367e-16\\
139	5.74737336103367e-16\\
140	5.74737336103367e-16\\
141	5.74737336103367e-16\\
142	5.74737336103367e-16\\
143	5.74737336103367e-16\\
144	5.74737336103367e-16\\
145	5.74737336103367e-16\\
146	5.74737336103367e-16\\
147	5.74737336103367e-16\\
148	5.74737336103367e-16\\
149	5.74737336103367e-16\\
150	5.74737336103367e-16\\
151	5.74737336103367e-16\\
152	5.74737336103367e-16\\
153	5.74737336103367e-16\\
154	5.74737336103367e-16\\
155	5.74737336103367e-16\\
156	5.74737336103367e-16\\
157	5.74737336103367e-16\\
158	5.74737336103367e-16\\
159	5.74737336103367e-16\\
160	5.74737336103367e-16\\
161	5.74737336103367e-16\\
162	5.74737336103367e-16\\
163	5.74737336103367e-16\\
164	5.74737336103367e-16\\
165	5.74737336103367e-16\\
166	5.74737336103367e-16\\
167	5.74737336103367e-16\\
168	5.74737336103367e-16\\
169	5.74737336103367e-16\\
170	5.74737336103367e-16\\
171	5.74737336103367e-16\\
172	5.74737336103367e-16\\
173	5.74737336103367e-16\\
174	5.74737336103367e-16\\
175	5.74737336103367e-16\\
176	5.74737336103367e-16\\
177	5.74737336103367e-16\\
178	5.74737336103367e-16\\
179	5.74737336103367e-16\\
180	5.74737336103367e-16\\
181	5.74737336103367e-16\\
182	5.74737336103367e-16\\
183	5.74737336103367e-16\\
184	5.74737336103367e-16\\
185	5.74737336103367e-16\\
186	5.74737336103367e-16\\
187	5.74737336103367e-16\\
188	5.74737336103367e-16\\
189	5.74737336103367e-16\\
190	5.74737336103367e-16\\
191	5.74737336103367e-16\\
192	5.74737336103367e-16\\
193	5.74737336103367e-16\\
194	5.74737336103367e-16\\
195	5.74737336103367e-16\\
196	5.74737336103367e-16\\
197	5.74737336103367e-16\\
198	5.74737336103367e-16\\
199	5.74737336103367e-16\\
200	5.74737336103367e-16\\
201	5.74737336103367e-16\\
202	5.74737336103367e-16\\
203	5.74737336103367e-16\\
204	5.74737336103367e-16\\
205	5.74737336103367e-16\\
206	5.74737336103367e-16\\
207	5.74737336103367e-16\\
208	5.74737336103367e-16\\
209	5.74737336103367e-16\\
210	5.74737336103367e-16\\
211	5.74737336103367e-16\\
212	5.74737336103367e-16\\
213	5.74737336103367e-16\\
214	5.74737336103367e-16\\
215	5.74737336103367e-16\\
216	5.74737336103367e-16\\
217	5.74737336103367e-16\\
218	5.74737336103367e-16\\
219	5.74737336103367e-16\\
220	5.74737336103367e-16\\
221	5.74737336103367e-16\\
222	5.74737336103367e-16\\
223	5.74737336103367e-16\\
224	5.74737336103367e-16\\
225	5.74737336103367e-16\\
226	5.74737336103367e-16\\
227	5.74737336103367e-16\\
228	5.74737336103367e-16\\
229	5.74737336103367e-16\\
230	5.74737336103367e-16\\
231	5.74737336103367e-16\\
232	5.74737336103367e-16\\
233	5.74737336103367e-16\\
234	5.74737336103367e-16\\
235	5.74737336103367e-16\\
236	5.74737336103367e-16\\
237	5.74737336103367e-16\\
238	5.74737336103367e-16\\
239	5.74737336103367e-16\\
240	5.74737336103367e-16\\
241	5.74737336103367e-16\\
242	5.74737336103367e-16\\
243	5.74737336103367e-16\\
244	5.74737336103367e-16\\
245	5.74737336103367e-16\\
246	5.74737336103367e-16\\
247	5.74737336103367e-16\\
248	5.74737336103367e-16\\
249	5.74737336103367e-16\\
250	5.74737336103367e-16\\
251	5.74737336103367e-16\\
252	5.74737336103367e-16\\
253	5.74737336103367e-16\\
254	5.74737336103367e-16\\
255	5.74737336103367e-16\\
256	5.74737336103367e-16\\
257	5.74737336103367e-16\\
258	5.74737336103367e-16\\
259	5.74737336103367e-16\\
260	5.74737336103367e-16\\
261	5.74737336103367e-16\\
262	5.74737336103367e-16\\
263	5.74737336103367e-16\\
264	5.74737336103367e-16\\
265	5.74737336103367e-16\\
266	5.74737336103367e-16\\
267	5.74737336103367e-16\\
268	5.74737336103367e-16\\
269	5.74737336103367e-16\\
270	5.74737336103367e-16\\
271	5.74737336103367e-16\\
272	5.74737336103367e-16\\
273	5.74737336103367e-16\\
274	5.74737336103367e-16\\
275	5.74737336103367e-16\\
276	5.74737336103367e-16\\
277	5.74737336103367e-16\\
278	5.74737336103367e-16\\
279	5.74737336103367e-16\\
280	5.74737336103367e-16\\
281	5.74737336103367e-16\\
282	5.74737336103367e-16\\
283	5.74737336103367e-16\\
284	5.74737336103367e-16\\
285	5.74737336103367e-16\\
286	5.74737336103367e-16\\
287	5.74737336103367e-16\\
288	5.74737336103367e-16\\
289	5.74737336103367e-16\\
290	5.74737336103367e-16\\
291	5.74737336103367e-16\\
292	5.74737336103367e-16\\
293	5.74737336103367e-16\\
294	5.74737336103367e-16\\
295	5.74737336103367e-16\\
296	5.74737336103367e-16\\
297	5.74737336103367e-16\\
298	5.74737336103367e-16\\
299	5.74737336103367e-16\\
300	5.74737336103367e-16\\
301	5.74737336103367e-16\\
302	5.74737336103367e-16\\
303	5.74737336103367e-16\\
304	5.74737336103367e-16\\
305	5.74737336103367e-16\\
306	5.74737336103367e-16\\
307	5.74737336103367e-16\\
308	5.74737336103367e-16\\
309	5.74737336103367e-16\\
310	5.74737336103367e-16\\
311	5.74737336103367e-16\\
312	5.74737336103367e-16\\
313	5.74737336103367e-16\\
314	5.74737336103367e-16\\
315	5.74737336103367e-16\\
316	5.74737336103367e-16\\
317	5.74737336103367e-16\\
318	5.74737336103367e-16\\
319	5.74737336103367e-16\\
320	5.74737336103367e-16\\
321	5.74737336103367e-16\\
322	5.74737336103367e-16\\
323	5.74737336103367e-16\\
324	5.74737336103367e-16\\
325	5.74737336103367e-16\\
326	5.74737336103367e-16\\
327	5.74737336103367e-16\\
328	5.74737336103367e-16\\
329	5.74737336103367e-16\\
330	5.74737336103367e-16\\
331	5.74737336103367e-16\\
332	5.74737336103367e-16\\
333	5.74737336103367e-16\\
334	5.74737336103367e-16\\
335	5.74737336103367e-16\\
336	5.74737336103367e-16\\
337	5.74737336103367e-16\\
338	5.74737336103367e-16\\
339	5.74737336103367e-16\\
340	5.74737336103367e-16\\
341	5.74737336103367e-16\\
342	5.74737336103367e-16\\
343	5.74737336103367e-16\\
344	5.74737336103367e-16\\
345	5.74737336103367e-16\\
346	5.74737336103367e-16\\
347	5.74737336103367e-16\\
348	5.74737336103367e-16\\
349	5.74737336103367e-16\\
350	5.74737336103367e-16\\
351	5.74737336103367e-16\\
352	5.74737336103367e-16\\
353	5.74737336103367e-16\\
354	5.74737336103367e-16\\
355	5.74737336103367e-16\\
356	5.74737336103367e-16\\
357	5.74737336103367e-16\\
358	5.74737336103367e-16\\
359	5.74737336103367e-16\\
360	5.74737336103367e-16\\
361	5.74737336103367e-16\\
362	5.74737336103367e-16\\
363	5.74737336103367e-16\\
364	5.74737336103367e-16\\
365	5.74737336103367e-16\\
366	5.74737336103367e-16\\
367	5.74737336103367e-16\\
368	5.74737336103367e-16\\
369	5.74737336103367e-16\\
370	5.74737336103367e-16\\
371	5.74737336103367e-16\\
372	5.74737336103367e-16\\
373	5.74737336103367e-16\\
374	5.74737336103367e-16\\
375	5.74737336103367e-16\\
376	5.74737336103367e-16\\
377	5.74737336103367e-16\\
378	5.74737336103367e-16\\
379	5.74737336103367e-16\\
380	5.74737336103367e-16\\
381	5.74737336103367e-16\\
382	5.74737336103367e-16\\
383	5.74737336103367e-16\\
384	5.74737336103367e-16\\
385	5.74737336103367e-16\\
386	5.74737336103367e-16\\
387	5.74737336103367e-16\\
388	5.74737336103367e-16\\
389	5.74737336103367e-16\\
390	5.74737336103367e-16\\
391	5.74737336103367e-16\\
392	5.74737336103367e-16\\
393	5.74737336103367e-16\\
394	5.74737336103367e-16\\
395	5.74737336103367e-16\\
396	5.74737336103367e-16\\
397	5.74737336103367e-16\\
398	4.19569295542163e-16\\
399	2.6485342705677e-16\\
400	1.11102313384135e-16\\
};
\addlegendentry{\texttt{OS}}

\addplot [color=mycolor2, line width=1.5pt]
  table[row sep=crcr]{%
1	0.20104188883882\\
2	0.0768090190854893\\
3	0.0334118497187814\\
4	0.0185188472416509\\
5	0.0142263640646382\\
6	0.00810977240772585\\
7	0.0057014224232884\\
8	0.00375083246556658\\
9	0.00372171693274864\\
10	0.00211629891434464\\
11	0.00171462341992441\\
12	0.00161256218157528\\
13	0.000872332926550155\\
14	0.000762963300350404\\
15	0.000726006112594861\\
16	0.000690248493845876\\
17	0.00037661098402797\\
18	0.000329108091759387\\
19	0.000278171548240165\\
20	0.000230744377019114\\
21	0.00021341261514358\\
22	0.000165963110425159\\
23	0.000131745491883742\\
24	0.000103642978728293\\
25	8.82566720686996e-05\\
26	7.18752551404859e-05\\
27	6.94527449350742e-05\\
28	5.2536527909288e-05\\
29	4.8635743868839e-05\\
30	3.74891710713221e-05\\
31	3.53190274347578e-05\\
32	2.68839824015661e-05\\
33	2.19039236028896e-05\\
34	2.03505391879159e-05\\
35	1.66948495459209e-05\\
36	1.49843039092802e-05\\
37	1.26752919880038e-05\\
38	1.20243787579034e-05\\
39	9.02935407106381e-06\\
40	7.38645427505429e-06\\
41	6.62038830236343e-06\\
42	5.97328758810327e-06\\
43	5.18019752583128e-06\\
44	4.44736035026139e-06\\
45	3.98652802031584e-06\\
46	3.10988472921011e-06\\
47	2.69693783649443e-06\\
48	2.50028632899986e-06\\
49	2.05548049756455e-06\\
50	1.82801538085205e-06\\
51	1.54121580121871e-06\\
52	1.3791273582593e-06\\
53	1.21206725805878e-06\\
54	1.1694284053923e-06\\
55	8.50757316022162e-07\\
56	7.57906682388904e-07\\
57	6.3730520655355e-07\\
58	6.10021188557263e-07\\
59	4.90491714329934e-07\\
60	4.61285100635284e-07\\
61	3.97431911427087e-07\\
62	3.48609747124866e-07\\
63	3.31934960260492e-07\\
64	2.51751196509212e-07\\
65	2.32108004485793e-07\\
66	1.96899327021769e-07\\
67	1.83289467317167e-07\\
68	1.4887519748071e-07\\
69	1.27304889401261e-07\\
70	1.15164606660768e-07\\
71	1.03145269401743e-07\\
72	8.72068111923934e-08\\
73	8.21235403903338e-08\\
74	7.3673200351975e-08\\
75	6.66136212749968e-08\\
76	5.4123942726585e-08\\
77	4.76569320831402e-08\\
78	4.2618048395712e-08\\
79	3.73496478268215e-08\\
80	2.98025324735991e-08\\
81	2.59152343039842e-08\\
82	2.42037347891977e-08\\
83	2.00692188311975e-08\\
84	1.90561157147251e-08\\
85	1.83495320229344e-08\\
86	1.5902962219436e-08\\
87	1.42726316747919e-08\\
88	1.17178943830147e-08\\
89	1.08115567219219e-08\\
90	8.86350308537395e-09\\
91	7.93517824951507e-09\\
92	6.54433734229278e-09\\
93	6.32313977838604e-09\\
94	5.41536553041295e-09\\
95	4.27029398136919e-09\\
96	4.03639765859093e-09\\
97	3.7413226834562e-09\\
98	3.51054977558643e-09\\
99	3.20807541080663e-09\\
100	2.9732982910379e-09\\
101	2.85727289748971e-09\\
102	2.3360472038473e-09\\
103	2.18316564523591e-09\\
104	2.12188040066994e-09\\
105	2.02686370188837e-09\\
106	1.90299533552136e-09\\
107	1.84894338361567e-09\\
108	1.67250201547573e-09\\
109	1.63938201776455e-09\\
110	1.56627141039255e-09\\
111	1.52926008523184e-09\\
112	1.48154211406722e-09\\
113	1.46581673873228e-09\\
114	1.44157996881146e-09\\
115	1.42040418266337e-09\\
116	1.40058866391731e-09\\
117	1.37212483325219e-09\\
118	1.33231666377066e-09\\
119	1.31725448349308e-09\\
120	1.30765880048359e-09\\
121	1.28833516245634e-09\\
122	1.26815214597818e-09\\
123	1.25658312634089e-09\\
124	1.2402010505601e-09\\
125	1.21913375989666e-09\\
126	1.20292634549531e-09\\
127	1.1749400312602e-09\\
128	1.15556687745691e-09\\
129	1.13624708998892e-09\\
130	1.12475064729103e-09\\
131	1.1167240776338e-09\\
132	1.09867401606636e-09\\
133	1.08379467593984e-09\\
134	1.07534855973529e-09\\
135	1.06465176862185e-09\\
136	1.04168682126685e-09\\
137	1.03600887863277e-09\\
138	1.02737219795729e-09\\
139	1.00563279832173e-09\\
140	1.00193732349673e-09\\
141	9.86008362581207e-10\\
142	9.71280237410399e-10\\
143	9.66233974077888e-10\\
144	9.52219541405511e-10\\
145	9.44400185848041e-10\\
146	9.35808722185329e-10\\
147	9.30195827142544e-10\\
148	9.17222770915272e-10\\
149	9.0401671127063e-10\\
150	8.99747868276638e-10\\
151	8.88089974948398e-10\\
152	8.76006964106482e-10\\
153	8.71054889904154e-10\\
154	8.61984015530455e-10\\
155	8.55378743810339e-10\\
156	8.37045013213556e-10\\
157	8.22517328549049e-10\\
158	8.15672007621384e-10\\
159	7.99346863155301e-10\\
160	7.98480985846836e-10\\
161	7.94050670250794e-10\\
162	7.73774934699502e-10\\
163	7.67148124382544e-10\\
164	7.62677475538579e-10\\
165	7.5617798212204e-10\\
166	7.51100410530026e-10\\
167	7.42436017682841e-10\\
168	7.29145869991972e-10\\
169	7.25278245832915e-10\\
170	7.19315785782491e-10\\
171	7.18801594309519e-10\\
172	7.07872282945491e-10\\
173	7.0236202202683e-10\\
174	6.99169014875125e-10\\
175	6.90191774310536e-10\\
176	6.83524911698073e-10\\
177	6.76862659370989e-10\\
178	6.71626184141186e-10\\
179	6.65295605503049e-10\\
180	6.58902845391545e-10\\
181	6.51600376413702e-10\\
182	6.45088006742035e-10\\
183	6.34036431132796e-10\\
184	6.30828617880073e-10\\
185	6.25065249942589e-10\\
186	6.18717813364393e-10\\
187	6.10451738302791e-10\\
188	6.03495571228926e-10\\
189	6.01240295704977e-10\\
190	5.91612127553522e-10\\
191	5.89555397943668e-10\\
192	5.78675633030463e-10\\
193	5.69402308278541e-10\\
194	5.68610179969371e-10\\
195	5.61625978184987e-10\\
196	5.55331458857912e-10\\
197	5.52890569429866e-10\\
198	5.47329440023582e-10\\
199	5.37768903207523e-10\\
200	5.33252298478842e-10\\
201	5.32037188084482e-10\\
202	5.26295422457417e-10\\
203	5.21742061716922e-10\\
204	5.19542688521956e-10\\
205	5.06594208761663e-10\\
206	5.02070727918297e-10\\
207	4.98869682391951e-10\\
208	4.87959975827728e-10\\
209	4.83755407680217e-10\\
210	4.81198703177917e-10\\
211	4.77100884893893e-10\\
212	4.73533802720983e-10\\
213	4.71616025382548e-10\\
214	4.67574041861908e-10\\
215	4.61543760866632e-10\\
216	4.60255672574644e-10\\
217	4.49687496932433e-10\\
218	4.46154280410731e-10\\
219	4.38521522971337e-10\\
220	4.33615171546153e-10\\
221	4.29880437874149e-10\\
222	4.26412802974088e-10\\
223	4.23197361992865e-10\\
224	4.17175614116807e-10\\
225	4.15108652369047e-10\\
226	4.1334239063697e-10\\
227	4.08060976662702e-10\\
228	4.0314480679599e-10\\
229	4.00891932620819e-10\\
230	3.97900573654396e-10\\
231	3.90034510839632e-10\\
232	3.87013769971558e-10\\
233	3.84115312383495e-10\\
234	3.79594260379151e-10\\
235	3.73255995563044e-10\\
236	3.69279392216042e-10\\
237	3.66208808396593e-10\\
238	3.62637698416275e-10\\
239	3.59668110304789e-10\\
240	3.56027271092022e-10\\
241	3.50315241304895e-10\\
242	3.48671815989535e-10\\
243	3.45239360429505e-10\\
244	3.44945307263351e-10\\
245	3.42774667851684e-10\\
246	3.37485840096839e-10\\
247	3.3655891698018e-10\\
248	3.3151387967436e-10\\
249	3.26421399051624e-10\\
250	3.22127942241216e-10\\
251	3.18244153273004e-10\\
252	3.18059092657673e-10\\
253	3.11852285686498e-10\\
254	3.05814813656371e-10\\
255	3.02697899727528e-10\\
256	3.01406288401051e-10\\
257	2.98040237913374e-10\\
258	2.95336116093726e-10\\
259	2.94083907218897e-10\\
260	2.92330344484418e-10\\
261	2.88730666584358e-10\\
262	2.83738845617383e-10\\
263	2.82762083422205e-10\\
264	2.81359950561627e-10\\
265	2.72996647264419e-10\\
266	2.70513385643333e-10\\
267	2.66602573650299e-10\\
268	2.64822339372712e-10\\
269	2.63649676879951e-10\\
270	2.58783616401732e-10\\
271	2.56448908624368e-10\\
272	2.53406036812736e-10\\
273	2.50391641461331e-10\\
274	2.48911386094254e-10\\
275	2.48226169046739e-10\\
276	2.4361118402857e-10\\
277	2.43296691703833e-10\\
278	2.41430265417424e-10\\
279	2.38441916104023e-10\\
280	2.35664181337454e-10\\
281	2.31080508957493e-10\\
282	2.25993729317449e-10\\
283	2.25154822351067e-10\\
284	2.21690990026188e-10\\
285	2.19331232798804e-10\\
286	2.17964119915529e-10\\
287	2.15206188383511e-10\\
288	2.12263019253013e-10\\
289	2.10708192980561e-10\\
290	2.08568375937417e-10\\
291	2.06769734924716e-10\\
292	2.01974049044117e-10\\
293	2.00146218586075e-10\\
294	1.98933818038501e-10\\
295	1.94535181532389e-10\\
296	1.93107375910023e-10\\
297	1.91951597159998e-10\\
298	1.88828068686405e-10\\
299	1.88371335654677e-10\\
300	1.84161308151399e-10\\
301	1.82988093307733e-10\\
302	1.79814729591397e-10\\
303	1.76161007357323e-10\\
304	1.73800458935548e-10\\
305	1.71541705442649e-10\\
306	1.69334655603799e-10\\
307	1.68276030143529e-10\\
308	1.6639471669356e-10\\
309	1.66314837680002e-10\\
310	1.65181893538641e-10\\
311	1.61955077664233e-10\\
312	1.58513086573062e-10\\
313	1.56592694842498e-10\\
314	1.55567490574488e-10\\
315	1.51901922933543e-10\\
316	1.50871811483235e-10\\
317	1.46153175845459e-10\\
318	1.46033932352738e-10\\
319	1.45342870930159e-10\\
320	1.43387680307845e-10\\
321	1.40511879691201e-10\\
322	1.398593632301e-10\\
323	1.37820129603921e-10\\
324	1.36453595272618e-10\\
325	1.35053506325509e-10\\
326	1.32922844141006e-10\\
327	1.28811283499368e-10\\
328	1.27953676905881e-10\\
329	1.25659839427782e-10\\
330	1.25295344350981e-10\\
331	1.24498910856388e-10\\
332	1.20858968679528e-10\\
333	1.20268254112003e-10\\
334	1.17000963343512e-10\\
335	1.16575354404907e-10\\
336	1.14268366703036e-10\\
337	1.13796224374309e-10\\
338	1.11216071148469e-10\\
339	1.09770978768776e-10\\
340	1.09308103580039e-10\\
341	1.07591374829528e-10\\
342	1.05301754517733e-10\\
343	1.03873076452688e-10\\
344	1.01634387945749e-10\\
345	1.00422158938708e-10\\
346	9.84178770694751e-11\\
347	9.71690304654361e-11\\
348	9.40711850896563e-11\\
349	9.30565119550458e-11\\
350	9.12565972217772e-11\\
351	8.99053072404326e-11\\
352	8.8008326323408e-11\\
353	8.68943714958063e-11\\
354	8.53005692783321e-11\\
355	8.4162637202255e-11\\
356	8.08856267639847e-11\\
357	8.04389261982758e-11\\
358	7.98928851630448e-11\\
359	7.85702412030834e-11\\
360	7.6868336351923e-11\\
361	7.54945326034126e-11\\
362	7.28365524599319e-11\\
363	7.21552899490962e-11\\
364	7.08805907671931e-11\\
365	6.94322057438868e-11\\
366	6.70514711630969e-11\\
367	6.55642505287097e-11\\
368	6.48973689728501e-11\\
369	6.35604244145878e-11\\
370	6.24289197883142e-11\\
371	6.06771198708661e-11\\
372	5.8746718952956e-11\\
373	5.6762643979714e-11\\
374	5.53426735394293e-11\\
375	5.35554268530902e-11\\
376	5.06482523958184e-11\\
377	4.97775972241313e-11\\
378	4.83942266773942e-11\\
379	4.74530511915652e-11\\
380	4.60244932125217e-11\\
381	4.44567424746462e-11\\
382	4.36370528426327e-11\\
383	4.18446429138666e-11\\
384	3.89054472170665e-11\\
385	3.72015541183456e-11\\
386	3.68704137008255e-11\\
387	3.56451591908027e-11\\
388	3.37448158063557e-11\\
389	3.20563435809792e-11\\
390	3.09368152319026e-11\\
391	2.93542750509134e-11\\
392	2.76703193834481e-11\\
393	2.63348667626326e-11\\
394	2.3872421109093e-11\\
395	2.23238727945347e-11\\
396	1.64405259228381e-11\\
397	1.48853678519777e-11\\
398	1.37864331385219e-11\\
399	1.02272370422199e-11\\
400	5.39364226219189e-12\\
};
\addlegendentry{\texttt{BT}}

\addplot [color=black, dashed, line width=1.5pt, forget plot]
  table[row sep=crcr]{%
25	1e-20\\
25	6.76440006564327\\
};
\end{axis}
\end{tikzpicture}%

%% file: stochastic_os_vs_bt_r25_nt256_NS1e5.tex
%
%
\definecolor{mycolor1}{rgb}{0.00000,0.44700,0.74100}%
\definecolor{mycolor2}{rgb}{0.85000,0.32500,0.09800}%
\definecolor{mycolor3}{rgb}{0.92900,0.69400,0.12500}%
\begin{tikzpicture}

\begin{axis}[%
width=0.951\wex,
height=\hex,
at={(0\wex,0\hex)},
scale only axis,
xmin=0,
xmax=1,
xlabel style={font=\color{white!15!black}},
xlabel={time (t)},
ymode=log,
ymin=1e-08,
ymax=0.009,
yminorticks=true,
ylabel style={font=\color{white!15!black}},
ylabel={Absolute mean value error},
axis background/.style={fill=white},
title style={font=\bfseries},
title={Time domain expected values simulation},
legend style={legend cell align=left, align=left, draw=white!15!black},
legend pos = south east,
xmajorgrids,
ymajorgrids,
yminorgrids 
]
\addplot [color=mycolor1, line width=1.5pt]
  table[row sep=crcr]{%
0	0\\
0.00390625	0\\
0.0078125	6.33612660787501e-09\\
0.01171875	4.63101493101759e-08\\
0.015625	1.77725547991382e-07\\
0.01953125	4.41207906710779e-07\\
0.0234375	7.70033491103206e-07\\
0.02734375	1.12174324282577e-06\\
0.03125	1.68313113525848e-06\\
0.03515625	2.72709868991346e-06\\
0.0390625	4.20345927771478e-06\\
0.04296875	5.92677175171824e-06\\
0.046875	7.70760082576846e-06\\
0.05078125	9.39963602740164e-06\\
0.0546875	1.09880965664716e-05\\
0.05859375	1.25160558397465e-05\\
0.0625	1.41607974253302e-05\\
0.06640625	1.61660596112724e-05\\
0.0703125	1.86578363805916e-05\\
0.07421875	2.16827390062754e-05\\
0.078125	2.51467606966544e-05\\
0.08203125	2.89443279756451e-05\\
0.0859375	3.29114397799107e-05\\
0.08984375	3.69367361951756e-05\\
0.09375	4.09826280680015e-05\\
0.09765625	4.49768966332192e-05\\
0.1015625	4.87839219599812e-05\\
0.10546875	5.23587963645259e-05\\
0.109375	5.57276774117955e-05\\
0.11328125	5.88239916530693e-05\\
0.1171875	6.17459022552233e-05\\
0.12109375	6.44285861881992e-05\\
0.125	6.68822700023525e-05\\
0.12890625	6.91847896205807e-05\\
0.1328125	7.13865275710628e-05\\
0.13671875	7.36108068874644e-05\\
0.140625	7.58888078414951e-05\\
0.14453125	7.81892501331843e-05\\
0.1484375	8.0769075527186e-05\\
0.15234375	8.36529319248945e-05\\
0.15625	8.68530696613547e-05\\
0.16015625	9.04585208439703e-05\\
0.1640625	9.46322049067201e-05\\
0.16796875	9.96190692618359e-05\\
0.171875	0.000105341948636939\\
0.17578125	0.000111605917886821\\
0.1796875	0.000118221238091622\\
0.18359375	0.00012531888055661\\
0.1875	0.000132915594705472\\
0.19140625	0.000140816876849782\\
0.1953125	0.000149093896448219\\
0.19921875	0.000157532514236608\\
0.203125	0.000166221374002578\\
0.20703125	0.000174954396040973\\
0.2109375	0.000183746041071824\\
0.21484375	0.000192778306084246\\
0.21875	0.000201784837127355\\
0.22265625	0.000210627493079507\\
0.2265625	0.000219562151190278\\
0.23046875	0.000228454745839964\\
0.234375	0.000237287609266822\\
0.23828125	0.000246020884428551\\
0.2421875	0.000254383279987543\\
0.24609375	0.000262692846626329\\
0.25	0.000270796130765409\\
0.25390625	0.000278591005708376\\
0.2578125	0.000286063134145938\\
0.26171875	0.000293120412036254\\
0.265625	0.000299740843074376\\
0.26953125	0.000306172719077229\\
0.2734375	0.00031231761871124\\
0.27734375	0.000318092516286055\\
0.28125	0.000323484532017937\\
0.28515625	0.000328309490214323\\
0.2890625	0.000332774272958324\\
0.29296875	0.000336935078914358\\
0.296875	0.00034057276966194\\
0.30078125	0.000343996982531037\\
0.3046875	0.000347068939505657\\
0.30859375	0.000349701432951079\\
0.3125	0.000351892330043583\\
0.31640625	0.000353600545463306\\
0.3203125	0.000354834218501637\\
0.32421875	0.00035565017298711\\
0.328125	0.000356073193854329\\
0.33203125	0.000356136341933794\\
0.3359375	0.000355586841888916\\
0.33984375	0.000354494001790296\\
0.34375	0.000352997924990589\\
0.34765625	0.000351018246784061\\
0.3515625	0.000348434377098017\\
0.35546875	0.00034543397673359\\
0.359375	0.000342117696328457\\
0.36328125	0.000338354185980669\\
0.3671875	0.000334341680643628\\
0.37109375	0.000329641704458599\\
0.375	0.000324652472855405\\
0.37890625	0.000319252818488547\\
0.3828125	0.000313487917388575\\
0.38671875	0.000307297992528781\\
0.390625	0.000300863894160799\\
0.39453125	0.00029437293602321\\
0.3984375	0.00028764937039054\\
0.40234375	0.000280657964476554\\
0.40625	0.000273414836927357\\
0.41015625	0.000265991212496403\\
0.4140625	0.000258453291895449\\
0.41796875	0.000250608915124877\\
0.421875	0.000242752966723207\\
0.42578125	0.000235147599462311\\
0.4296875	0.000227243020397391\\
0.43359375	0.000219210056883606\\
0.4375	0.000209354129989201\\
0.44140625	0.000199116139776233\\
0.4453125	0.000189006649882501\\
0.44921875	0.000179214403910917\\
0.453125	0.000169672530225021\\
0.45703125	0.000160406916920278\\
0.4609375	0.000151609084858079\\
0.46484375	0.000143245701032206\\
0.46875	0.000135423652341218\\
0.47265625	0.000128104646246936\\
0.4765625	0.000121536427713882\\
0.48046875	0.000115655245483882\\
0.484375	0.000110741263064471\\
0.48828125	0.000106888892064178\\
0.4921875	0.000104150958418929\\
0.49609375	0.00010246415971047\\
0.5	0.000102107220109324\\
0.50390625	0.000103057578699129\\
0.5078125	0.000105487212371648\\
0.51171875	0.000109275535190131\\
0.515625	0.000114414930195822\\
0.51953125	0.000120499741176566\\
0.5234375	0.000127385728184325\\
0.52734375	0.000134901977592194\\
0.53125	0.000142791327662084\\
0.53515625	0.000150935833666484\\
0.5390625	0.000159150014686805\\
0.54296875	0.000167449530713521\\
0.546875	0.000175861927289007\\
0.55078125	0.000184110950233154\\
0.5546875	0.000192244549135788\\
0.55859375	0.000200188724125525\\
0.5625	0.000207926810227308\\
0.56640625	0.000215449677862276\\
0.5703125	0.00022278134508708\\
0.57421875	0.000229845578126041\\
0.578125	0.00023668026801695\\
0.58203125	0.000243183930944851\\
0.5859375	0.000249410756415357\\
0.58984375	0.00025529345271331\\
0.59375	0.00026096060407499\\
0.59765625	0.00026642859165011\\
0.6015625	0.000271285603447758\\
0.60546875	0.000275711661838147\\
0.609375	0.000279776009991248\\
0.61328125	0.000283535526010042\\
0.6171875	0.000286733522039737\\
0.62109375	0.000289555700432467\\
0.625	0.000291803444435531\\
0.62890625	0.000293816011411242\\
0.6328125	0.000295501097929308\\
0.63671875	0.00029663830054283\\
0.640625	0.000297277472411113\\
0.64453125	0.000297550476482922\\
0.6484375	0.000297691690965754\\
0.65234375	0.000297395228598469\\
0.65625	0.000296502597832373\\
0.66015625	0.00029534778401416\\
0.6640625	0.00029369701780469\\
0.66796875	0.000291487522155385\\
0.671875	0.000288965297700462\\
0.67578125	0.00028618755152729\\
0.6796875	0.000282858589723744\\
0.68359375	0.000279189247454795\\
0.6875	0.000275210304444511\\
0.69140625	0.000271045674803818\\
0.6953125	0.000266564213742602\\
0.69921875	0.000261784862923555\\
0.703125	0.000256580346321533\\
0.70703125	0.000250975546631871\\
0.7109375	0.000245158837803082\\
0.71484375	0.000238959604771853\\
0.71875	0.000232538535938871\\
0.72265625	0.000225834238972126\\
0.7265625	0.000218748109970491\\
0.73046875	0.000211408828817514\\
0.734375	0.000203852185485709\\
0.73828125	0.000196244927368648\\
0.7421875	0.000188598708330317\\
0.74609375	0.000180851631288384\\
0.75	0.000173066247013039\\
0.75390625	0.000165335631934819\\
0.7578125	0.000157544507151121\\
0.76171875	0.000149795189209068\\
0.765625	0.000142015625014747\\
0.76953125	0.00013453276751299\\
0.7734375	0.000127298901271129\\
0.77734375	0.000120342389208908\\
0.78125	0.000113779352758274\\
0.78515625	0.000107602080310382\\
0.7890625	0.00010194633621498\\
0.79296875	9.69236196138668e-05\\
0.796875	9.26697437382578e-05\\
0.80078125	8.91390955594006e-05\\
0.8046875	8.64658653579665e-05\\
0.80859375	8.48193188756459e-05\\
0.8125	8.43299606967983e-05\\
0.81640625	8.4916791320559e-05\\
0.8203125	8.67051018099577e-05\\
0.82421875	8.9680387928354e-05\\
0.828125	9.3736320071563e-05\\
0.83203125	9.87939156429126e-05\\
0.8359375	0.000104543315752379\\
0.83984375	0.000110863611756411\\
0.84375	0.00011753244445821\\
0.84765625	0.000124519176911262\\
0.8515625	0.000131587545462126\\
0.85546875	0.00013876472411482\\
0.859375	0.000145979293664397\\
0.86328125	0.000153052391949032\\
0.8671875	0.000160070698593654\\
0.87109375	0.000166999028489967\\
0.875	0.000173823127675231\\
0.87890625	0.000180540758990452\\
0.8828125	0.000187089658725486\\
0.88671875	0.000193442026372355\\
0.890625	0.000199408816641703\\
0.89453125	0.00020505608083603\\
0.8984375	0.000210624742341765\\
0.90234375	0.000215830749150097\\
0.90625	0.000220728676474835\\
0.91015625	0.000225383446379273\\
0.9140625	0.000229720876620941\\
0.91796875	0.000233878209420594\\
0.921875	0.000237527466887378\\
0.92578125	0.00024086903816286\\
0.9296875	0.00024392219077043\\
0.93359375	0.000246783632608417\\
0.9375	0.000249000530342076\\
0.94140625	0.000250839606576116\\
0.9453125	0.000252547120360335\\
0.94921875	0.000253774318677812\\
0.953125	0.000254395140234441\\
0.95703125	0.000254682357402357\\
0.9609375	0.000254814568299271\\
0.96484375	0.000254618102796482\\
0.96875	0.00025409594548844\\
0.97265625	0.000253148701370136\\
0.9765625	0.000251615828373965\\
0.98046875	0.000250059444393488\\
0.984375	0.000248146539251238\\
0.98828125	0.000245917794446086\\
0.9921875	0.00024331338116578\\
0.99609375	0.000240450239378557\\
1	0.000237180925607609\\
};
\addlegendentry{\texttt{OS}}

\addplot [color=mycolor2, line width=1.5pt]
  table[row sep=crcr]{%
0	0\\
0.00390625	0\\
0.0078125	3.3462118329318e-09\\
0.01171875	3.98899584309525e-08\\
0.015625	1.271882394901e-07\\
0.01953125	3.00670891513877e-07\\
0.0234375	5.51879309211068e-07\\
0.02734375	8.23142554703089e-07\\
0.03125	1.16859052792536e-06\\
0.03515625	1.75641850206756e-06\\
0.0390625	2.63784953061893e-06\\
0.04296875	3.70809038178176e-06\\
0.046875	4.84882264797916e-06\\
0.05078125	5.94795500441802e-06\\
0.0546875	6.97550723646515e-06\\
0.05859375	7.96244206785536e-06\\
0.0625	8.97485058762768e-06\\
0.06640625	1.01791740144081e-05\\
0.0703125	1.16578981762731e-05\\
0.07421875	1.34711911701868e-05\\
0.078125	1.55654985365582e-05\\
0.08203125	1.7885790752302e-05\\
0.0859375	2.03557952043055e-05\\
0.08984375	2.28659710727182e-05\\
0.09375	2.5369745749927e-05\\
0.09765625	2.78623150224234e-05\\
0.1015625	3.02623283389891e-05\\
0.10546875	3.24839679531008e-05\\
0.109375	3.46026168080097e-05\\
0.11328125	3.65099943728485e-05\\
0.1171875	3.8291217762671e-05\\
0.12109375	3.99572241688304e-05\\
0.125	4.14412244123991e-05\\
0.12890625	4.27780989904892e-05\\
0.1328125	4.40521192803171e-05\\
0.13671875	4.52597306804349e-05\\
0.140625	4.65294736860666e-05\\
0.14453125	4.7746505383504e-05\\
0.1484375	4.90664517425843e-05\\
0.15234375	5.05931355787627e-05\\
0.15625	5.22506647857341e-05\\
0.16015625	5.41010420984383e-05\\
0.1640625	5.61885903724387e-05\\
0.16796875	5.8765787576305e-05\\
0.171875	6.18189607178495e-05\\
0.17578125	6.52292117824587e-05\\
0.1796875	6.8819032026425e-05\\
0.18359375	7.27213131969623e-05\\
0.1875	7.69035590818629e-05\\
0.19140625	8.13161101481879e-05\\
0.1953125	8.58968591233615e-05\\
0.19921875	9.06146130906853e-05\\
0.203125	9.54988424248078e-05\\
0.20703125	0.000100484352320411\\
0.2109375	0.000105464811258812\\
0.21484375	0.000110506464814352\\
0.21875	0.000115654808497701\\
0.22265625	0.000120667691763577\\
0.2265625	0.000125739879945937\\
0.23046875	0.00013087065666353\\
0.234375	0.000135845402519643\\
0.23828125	0.000140842124068099\\
0.2421875	0.000145693293300704\\
0.24609375	0.000150451089178719\\
0.25	0.000155067760870299\\
0.25390625	0.000159531894311226\\
0.2578125	0.00016383503356834\\
0.26171875	0.000167935228824639\\
0.265625	0.000171788139846286\\
0.26953125	0.000175496938956558\\
0.2734375	0.000178977337698735\\
0.27734375	0.0001822448536882\\
0.28125	0.000185429607574718\\
0.28515625	0.000188247528840364\\
0.2890625	0.0001907901871666\\
0.29296875	0.000193273209811772\\
0.296875	0.000195369784291224\\
0.30078125	0.000197298453021971\\
0.3046875	0.000199108820801857\\
0.30859375	0.000200667917068822\\
0.3125	0.000201909389903256\\
0.31640625	0.000202991506426021\\
0.3203125	0.000203739696496219\\
0.32421875	0.000204300442007798\\
0.328125	0.00020456245452261\\
0.33203125	0.000204624786078144\\
0.3359375	0.0002043997627174\\
0.33984375	0.000203845647365203\\
0.34375	0.000203031063464506\\
0.34765625	0.000201933686507808\\
0.3515625	0.000200569211373053\\
0.35546875	0.000198850522736143\\
0.359375	0.000196973326257617\\
0.36328125	0.000194739045789046\\
0.3671875	0.000192574351318585\\
0.37109375	0.000189946564083324\\
0.375	0.000187065254346575\\
0.37890625	0.000183980568299159\\
0.3828125	0.000180747384948478\\
0.38671875	0.000177250625601119\\
0.390625	0.000173630508532273\\
0.39453125	0.000169844916463305\\
0.3984375	0.000166002678496697\\
0.40234375	0.000162004870966474\\
0.40625	0.000157852995250256\\
0.41015625	0.000153537940941052\\
0.4140625	0.000149188344973062\\
0.41796875	0.000144602087733425\\
0.421875	0.000140026371188463\\
0.42578125	0.00013556022480549\\
0.4296875	0.000130894210034799\\
0.43359375	0.000126216003784134\\
0.4375	0.000121066242841627\\
0.44140625	0.000115486952854551\\
0.4453125	0.000109818138426175\\
0.44921875	0.000104264079522092\\
0.453125	9.88545434110666e-05\\
0.45703125	9.35172897236903e-05\\
0.4609375	8.84757570362012e-05\\
0.46484375	8.36645571654768e-05\\
0.46875	7.91130429287463e-05\\
0.47265625	7.48655063905092e-05\\
0.4765625	7.0991449891372e-05\\
0.48046875	6.75937882624844e-05\\
0.484375	6.46498331487137e-05\\
0.48828125	6.22915922499712e-05\\
0.4921875	6.05781948511941e-05\\
0.49609375	5.94453643894082e-05\\
0.5	5.90759036815471e-05\\
0.50390625	5.94325762323816e-05\\
0.5078125	6.0600710696327e-05\\
0.51171875	6.26310275659789e-05\\
0.515625	6.54606395171462e-05\\
0.51953125	6.89195598668279e-05\\
0.5234375	7.28083030717416e-05\\
0.52734375	7.70792473591685e-05\\
0.53125	8.15736473602035e-05\\
0.53515625	8.62576720289216e-05\\
0.5390625	9.10181734942306e-05\\
0.54296875	9.57470146002733e-05\\
0.546875	0.000100620464221774\\
0.55078125	0.000105336380285155\\
0.5546875	0.000110043386168703\\
0.55859375	0.000114623122137712\\
0.5625	0.000119127833159776\\
0.56640625	0.000123458417269009\\
0.5703125	0.0001276802565235\\
0.57421875	0.000131693970902888\\
0.578125	0.000135580720506343\\
0.58203125	0.000139358380113605\\
0.5859375	0.000142914294848333\\
0.58984375	0.000146305101656696\\
0.59375	0.000149524836107424\\
0.59765625	0.000152725444070384\\
0.6015625	0.000155561417984687\\
0.60546875	0.0001581674397336\\
0.609375	0.000160542854842489\\
0.61328125	0.000162766348698225\\
0.6171875	0.00016458378536867\\
0.62109375	0.000166281065864954\\
0.625	0.00016762594085081\\
0.62890625	0.000168761619213715\\
0.6328125	0.000169745861174691\\
0.63671875	0.000170475906301037\\
0.640625	0.000170883703899239\\
0.64453125	0.000171015412109848\\
0.6484375	0.000171043864023108\\
0.65234375	0.000170992148860669\\
0.65625	0.000170577510171319\\
0.66015625	0.000169934703144844\\
0.6640625	0.000169018535152819\\
0.66796875	0.000167850900186107\\
0.671875	0.000166398032470826\\
0.67578125	0.000164820547262906\\
0.6796875	0.000162977494410974\\
0.68359375	0.000160885290545289\\
0.6875	0.000158592659376559\\
0.69140625	0.00015626266566539\\
0.6953125	0.000153644487799281\\
0.69921875	0.000150993101860851\\
0.703125	0.000148056498320289\\
0.70703125	0.000144932374025559\\
0.7109375	0.0001416897108973\\
0.71484375	0.00013807319531147\\
0.71875	0.000134493468615471\\
0.72265625	0.000130711211275348\\
0.7265625	0.00012669292198276\\
0.73046875	0.000122492507516062\\
0.734375	0.000118189178109085\\
0.73828125	0.000113857410361872\\
0.7421875	0.000109483698010568\\
0.74609375	0.000105100081136795\\
0.75	0.000100651520271034\\
0.75390625	9.62518370835187e-05\\
0.7578125	9.18571210388396e-05\\
0.76171875	8.74702375135621e-05\\
0.765625	8.30312346665477e-05\\
0.76953125	7.8751633354668e-05\\
0.7734375	7.4589652422331e-05\\
0.77734375	7.06193698545243e-05\\
0.78125	6.68656915124845e-05\\
0.78515625	6.33060838578612e-05\\
0.7890625	6.00189610785959e-05\\
0.79296875	5.70352432424139e-05\\
0.796875	5.44977912814004e-05\\
0.80078125	5.23692174881133e-05\\
0.8046875	5.06800560591002e-05\\
0.80859375	4.9544381848163e-05\\
0.8125	4.90964209824521e-05\\
0.81640625	4.92621465237176e-05\\
0.8203125	5.00993536192889e-05\\
0.82421875	5.16501004814591e-05\\
0.828125	5.38445750205864e-05\\
0.83203125	5.66724323584917e-05\\
0.8359375	5.98919687470452e-05\\
0.83984375	6.34894507123701e-05\\
0.84375	6.73120388583736e-05\\
0.84765625	7.13078857092522e-05\\
0.8515625	7.54125978139544e-05\\
0.85546875	7.95056934912113e-05\\
0.859375	8.36144446002904e-05\\
0.86328125	8.76345716110779e-05\\
0.8671875	9.16823439305115e-05\\
0.87109375	9.56440772956225e-05\\
0.875	9.95694718476704e-05\\
0.87890625	0.000103482172661164\\
0.8828125	0.000107134002234306\\
0.88671875	0.000110818487858576\\
0.890625	0.000114318341504142\\
0.89453125	0.000117561551314436\\
0.8984375	0.000120764169853506\\
0.90234375	0.000123755500770969\\
0.90625	0.000126531332359902\\
0.91015625	0.000129269102806522\\
0.9140625	0.000131796892990513\\
0.91796875	0.000134299068979271\\
0.921875	0.00013639014923485\\
0.92578125	0.000138338539375354\\
0.9296875	0.000139994523680572\\
0.93359375	0.000141625595733649\\
0.9375	0.000142987962585336\\
0.94140625	0.000144066143445491\\
0.9453125	0.000145069222783521\\
0.94921875	0.000145871185590909\\
0.953125	0.000146353893447538\\
0.95703125	0.000146439988910954\\
0.9609375	0.000146506457833345\\
0.96484375	0.000146410128645201\\
0.96875	0.000146120598577373\\
0.97265625	0.000145684244917292\\
0.9765625	0.00014479464338324\\
0.98046875	0.000143930020530283\\
0.984375	0.000142804570574712\\
0.98828125	0.000141624894240615\\
0.9921875	0.00014017928458673\\
0.99609375	0.000138498549358911\\
1	0.000136741697661798\\
};
\addlegendentry{\texttt{BT}}

\addplot [color=black, dashed, line width=1.5pt]
  table[row sep=crcr]{%
0	0.00302418301324739\\
0.00390625	0.00302418301324739\\
0.0078125	0.00302418301324739\\
0.01171875	0.00302418301324739\\
0.015625	0.00302418301324739\\
0.01953125	0.00302418301324739\\
0.0234375	0.00302418301324739\\
0.02734375	0.00302418301324739\\
0.03125	0.00302418301324739\\
0.03515625	0.00302418301324739\\
0.0390625	0.00302418301324739\\
0.04296875	0.00302418301324739\\
0.046875	0.00302418301324739\\
0.05078125	0.00302418301324739\\
0.0546875	0.00302418301324739\\
0.05859375	0.00302418301324739\\
0.0625	0.00302418301324739\\
0.06640625	0.00302418301324739\\
0.0703125	0.00302418301324739\\
0.07421875	0.00302418301324739\\
0.078125	0.00302418301324739\\
0.08203125	0.00302418301324739\\
0.0859375	0.00302418301324739\\
0.08984375	0.00302418301324739\\
0.09375	0.00302418301324739\\
0.09765625	0.00302418301324739\\
0.1015625	0.00302418301324739\\
0.10546875	0.00302418301324739\\
0.109375	0.00302418301324739\\
0.11328125	0.00302418301324739\\
0.1171875	0.00302418301324739\\
0.12109375	0.00302418301324739\\
0.125	0.00302418301324739\\
0.12890625	0.00302418301324739\\
0.1328125	0.00302418301324739\\
0.13671875	0.00302418301324739\\
0.140625	0.00302418301324739\\
0.14453125	0.00302418301324739\\
0.1484375	0.00302418301324739\\
0.15234375	0.00302418301324739\\
0.15625	0.00302418301324739\\
0.16015625	0.00302418301324739\\
0.1640625	0.00302418301324739\\
0.16796875	0.00302418301324739\\
0.171875	0.00302418301324739\\
0.17578125	0.00302418301324739\\
0.1796875	0.00302418301324739\\
0.18359375	0.00302418301324739\\
0.1875	0.00302418301324739\\
0.19140625	0.00302418301324739\\
0.1953125	0.00302418301324739\\
0.19921875	0.00302418301324739\\
0.203125	0.00302418301324739\\
0.20703125	0.00302418301324739\\
0.2109375	0.00302418301324739\\
0.21484375	0.00302418301324739\\
0.21875	0.00302418301324739\\
0.22265625	0.00302418301324739\\
0.2265625	0.00302418301324739\\
0.23046875	0.00302418301324739\\
0.234375	0.00302418301324739\\
0.23828125	0.00302418301324739\\
0.2421875	0.00302418301324739\\
0.24609375	0.00302418301324739\\
0.25	0.00302418301324739\\
0.25390625	0.00302418301324739\\
0.2578125	0.00302418301324739\\
0.26171875	0.00302418301324739\\
0.265625	0.00302418301324739\\
0.26953125	0.00302418301324739\\
0.2734375	0.00302418301324739\\
0.27734375	0.00302418301324739\\
0.28125	0.00302418301324739\\
0.28515625	0.00302418301324739\\
0.2890625	0.00302418301324739\\
0.29296875	0.00302418301324739\\
0.296875	0.00302418301324739\\
0.30078125	0.00302418301324739\\
0.3046875	0.00302418301324739\\
0.30859375	0.00302418301324739\\
0.3125	0.00302418301324739\\
0.31640625	0.00302418301324739\\
0.3203125	0.00302418301324739\\
0.32421875	0.00302418301324739\\
0.328125	0.00302418301324739\\
0.33203125	0.00302418301324739\\
0.3359375	0.00302418301324739\\
0.33984375	0.00302418301324739\\
0.34375	0.00302418301324739\\
0.34765625	0.00302418301324739\\
0.3515625	0.00302418301324739\\
0.35546875	0.00302418301324739\\
0.359375	0.00302418301324739\\
0.36328125	0.00302418301324739\\
0.3671875	0.00302418301324739\\
0.37109375	0.00302418301324739\\
0.375	0.00302418301324739\\
0.37890625	0.00302418301324739\\
0.3828125	0.00302418301324739\\
0.38671875	0.00302418301324739\\
0.390625	0.00302418301324739\\
0.39453125	0.00302418301324739\\
0.3984375	0.00302418301324739\\
0.40234375	0.00302418301324739\\
0.40625	0.00302418301324739\\
0.41015625	0.00302418301324739\\
0.4140625	0.00302418301324739\\
0.41796875	0.00302418301324739\\
0.421875	0.00302418301324739\\
0.42578125	0.00302418301324739\\
0.4296875	0.00302418301324739\\
0.43359375	0.00302418301324739\\
0.4375	0.00302418301324739\\
0.44140625	0.00302418301324739\\
0.4453125	0.00302418301324739\\
0.44921875	0.00302418301324739\\
0.453125	0.00302418301324739\\
0.45703125	0.00302418301324739\\
0.4609375	0.00302418301324739\\
0.46484375	0.00302418301324739\\
0.46875	0.00302418301324739\\
0.47265625	0.00302418301324739\\
0.4765625	0.00302418301324739\\
0.48046875	0.00302418301324739\\
0.484375	0.00302418301324739\\
0.48828125	0.00302418301324739\\
0.4921875	0.00302418301324739\\
0.49609375	0.00302418301324739\\
0.5	0.00302418301324739\\
0.50390625	0.00302418301324739\\
0.5078125	0.00302418301324739\\
0.51171875	0.00302418301324739\\
0.515625	0.00302418301324739\\
0.51953125	0.00302418301324739\\
0.5234375	0.00302418301324739\\
0.52734375	0.00302418301324739\\
0.53125	0.00302418301324739\\
0.53515625	0.00302418301324739\\
0.5390625	0.00302418301324739\\
0.54296875	0.00302418301324739\\
0.546875	0.00302418301324739\\
0.55078125	0.00302418301324739\\
0.5546875	0.00302418301324739\\
0.55859375	0.00302418301324739\\
0.5625	0.00302418301324739\\
0.56640625	0.00302418301324739\\
0.5703125	0.00302418301324739\\
0.57421875	0.00302418301324739\\
0.578125	0.00302418301324739\\
0.58203125	0.00302418301324739\\
0.5859375	0.00302418301324739\\
0.58984375	0.00302418301324739\\
0.59375	0.00302418301324739\\
0.59765625	0.00302418301324739\\
0.6015625	0.00302418301324739\\
0.60546875	0.00302418301324739\\
0.609375	0.00302418301324739\\
0.61328125	0.00302418301324739\\
0.6171875	0.00302418301324739\\
0.62109375	0.00302418301324739\\
0.625	0.00302418301324739\\
0.62890625	0.00302418301324739\\
0.6328125	0.00302418301324739\\
0.63671875	0.00302418301324739\\
0.640625	0.00302418301324739\\
0.64453125	0.00302418301324739\\
0.6484375	0.00302418301324739\\
0.65234375	0.00302418301324739\\
0.65625	0.00302418301324739\\
0.66015625	0.00302418301324739\\
0.6640625	0.00302418301324739\\
0.66796875	0.00302418301324739\\
0.671875	0.00302418301324739\\
0.67578125	0.00302418301324739\\
0.6796875	0.00302418301324739\\
0.68359375	0.00302418301324739\\
0.6875	0.00302418301324739\\
0.69140625	0.00302418301324739\\
0.6953125	0.00302418301324739\\
0.69921875	0.00302418301324739\\
0.703125	0.00302418301324739\\
0.70703125	0.00302418301324739\\
0.7109375	0.00302418301324739\\
0.71484375	0.00302418301324739\\
0.71875	0.00302418301324739\\
0.72265625	0.00302418301324739\\
0.7265625	0.00302418301324739\\
0.73046875	0.00302418301324739\\
0.734375	0.00302418301324739\\
0.73828125	0.00302418301324739\\
0.7421875	0.00302418301324739\\
0.74609375	0.00302418301324739\\
0.75	0.00302418301324739\\
0.75390625	0.00302418301324739\\
0.7578125	0.00302418301324739\\
0.76171875	0.00302418301324739\\
0.765625	0.00302418301324739\\
0.76953125	0.00302418301324739\\
0.7734375	0.00302418301324739\\
0.77734375	0.00302418301324739\\
0.78125	0.00302418301324739\\
0.78515625	0.00302418301324739\\
0.7890625	0.00302418301324739\\
0.79296875	0.00302418301324739\\
0.796875	0.00302418301324739\\
0.80078125	0.00302418301324739\\
0.8046875	0.00302418301324739\\
0.80859375	0.00302418301324739\\
0.8125	0.00302418301324739\\
0.81640625	0.00302418301324739\\
0.8203125	0.00302418301324739\\
0.82421875	0.00302418301324739\\
0.828125	0.00302418301324739\\
0.83203125	0.00302418301324739\\
0.8359375	0.00302418301324739\\
0.83984375	0.00302418301324739\\
0.84375	0.00302418301324739\\
0.84765625	0.00302418301324739\\
0.8515625	0.00302418301324739\\
0.85546875	0.00302418301324739\\
0.859375	0.00302418301324739\\
0.86328125	0.00302418301324739\\
0.8671875	0.00302418301324739\\
0.87109375	0.00302418301324739\\
0.875	0.00302418301324739\\
0.87890625	0.00302418301324739\\
0.8828125	0.00302418301324739\\
0.88671875	0.00302418301324739\\
0.890625	0.00302418301324739\\
0.89453125	0.00302418301324739\\
0.8984375	0.00302418301324739\\
0.90234375	0.00302418301324739\\
0.90625	0.00302418301324739\\
0.91015625	0.00302418301324739\\
0.9140625	0.00302418301324739\\
0.91796875	0.00302418301324739\\
0.921875	0.00302418301324739\\
0.92578125	0.00302418301324739\\
0.9296875	0.00302418301324739\\
0.93359375	0.00302418301324739\\
0.9375	0.00302418301324739\\
0.94140625	0.00302418301324739\\
0.9453125	0.00302418301324739\\
0.94921875	0.00302418301324739\\
0.953125	0.00302418301324739\\
0.95703125	0.00302418301324739\\
0.9609375	0.00302418301324739\\
0.96484375	0.00302418301324739\\
0.96875	0.00302418301324739\\
0.97265625	0.00302418301324739\\
0.9765625	0.00302418301324739\\
0.98046875	0.00302418301324739\\
0.984375	0.00302418301324739\\
0.98828125	0.00302418301324739\\
0.9921875	0.00302418301324739\\
0.99609375	0.00302418301324739\\
1	0.00302418301324739\\
};
\addlegendentry{Error bound \texttt{OS}}

\addplot [color=mycolor3, dotted, line width=1.5pt]
  table[row sep=crcr]{%
0	0.00311856032842183\\
0.00390625	0.00311856032842183\\
0.0078125	0.00311856032842183\\
0.01171875	0.00311856032842183\\
0.015625	0.00311856032842183\\
0.01953125	0.00311856032842183\\
0.0234375	0.00311856032842183\\
0.02734375	0.00311856032842183\\
0.03125	0.00311856032842183\\
0.03515625	0.00311856032842183\\
0.0390625	0.00311856032842183\\
0.04296875	0.00311856032842183\\
0.046875	0.00311856032842183\\
0.05078125	0.00311856032842183\\
0.0546875	0.00311856032842183\\
0.05859375	0.00311856032842183\\
0.0625	0.00311856032842183\\
0.06640625	0.00311856032842183\\
0.0703125	0.00311856032842183\\
0.07421875	0.00311856032842183\\
0.078125	0.00311856032842183\\
0.08203125	0.00311856032842183\\
0.0859375	0.00311856032842183\\
0.08984375	0.00311856032842183\\
0.09375	0.00311856032842183\\
0.09765625	0.00311856032842183\\
0.1015625	0.00311856032842183\\
0.10546875	0.00311856032842183\\
0.109375	0.00311856032842183\\
0.11328125	0.00311856032842183\\
0.1171875	0.00311856032842183\\
0.12109375	0.00311856032842183\\
0.125	0.00311856032842183\\
0.12890625	0.00311856032842183\\
0.1328125	0.00311856032842183\\
0.13671875	0.00311856032842183\\
0.140625	0.00311856032842183\\
0.14453125	0.00311856032842183\\
0.1484375	0.00311856032842183\\
0.15234375	0.00311856032842183\\
0.15625	0.00311856032842183\\
0.16015625	0.00311856032842183\\
0.1640625	0.00311856032842183\\
0.16796875	0.00311856032842183\\
0.171875	0.00311856032842183\\
0.17578125	0.00311856032842183\\
0.1796875	0.00311856032842183\\
0.18359375	0.00311856032842183\\
0.1875	0.00311856032842183\\
0.19140625	0.00311856032842183\\
0.1953125	0.00311856032842183\\
0.19921875	0.00311856032842183\\
0.203125	0.00311856032842183\\
0.20703125	0.00311856032842183\\
0.2109375	0.00311856032842183\\
0.21484375	0.00311856032842183\\
0.21875	0.00311856032842183\\
0.22265625	0.00311856032842183\\
0.2265625	0.00311856032842183\\
0.23046875	0.00311856032842183\\
0.234375	0.00311856032842183\\
0.23828125	0.00311856032842183\\
0.2421875	0.00311856032842183\\
0.24609375	0.00311856032842183\\
0.25	0.00311856032842183\\
0.25390625	0.00311856032842183\\
0.2578125	0.00311856032842183\\
0.26171875	0.00311856032842183\\
0.265625	0.00311856032842183\\
0.26953125	0.00311856032842183\\
0.2734375	0.00311856032842183\\
0.27734375	0.00311856032842183\\
0.28125	0.00311856032842183\\
0.28515625	0.00311856032842183\\
0.2890625	0.00311856032842183\\
0.29296875	0.00311856032842183\\
0.296875	0.00311856032842183\\
0.30078125	0.00311856032842183\\
0.3046875	0.00311856032842183\\
0.30859375	0.00311856032842183\\
0.3125	0.00311856032842183\\
0.31640625	0.00311856032842183\\
0.3203125	0.00311856032842183\\
0.32421875	0.00311856032842183\\
0.328125	0.00311856032842183\\
0.33203125	0.00311856032842183\\
0.3359375	0.00311856032842183\\
0.33984375	0.00311856032842183\\
0.34375	0.00311856032842183\\
0.34765625	0.00311856032842183\\
0.3515625	0.00311856032842183\\
0.35546875	0.00311856032842183\\
0.359375	0.00311856032842183\\
0.36328125	0.00311856032842183\\
0.3671875	0.00311856032842183\\
0.37109375	0.00311856032842183\\
0.375	0.00311856032842183\\
0.37890625	0.00311856032842183\\
0.3828125	0.00311856032842183\\
0.38671875	0.00311856032842183\\
0.390625	0.00311856032842183\\
0.39453125	0.00311856032842183\\
0.3984375	0.00311856032842183\\
0.40234375	0.00311856032842183\\
0.40625	0.00311856032842183\\
0.41015625	0.00311856032842183\\
0.4140625	0.00311856032842183\\
0.41796875	0.00311856032842183\\
0.421875	0.00311856032842183\\
0.42578125	0.00311856032842183\\
0.4296875	0.00311856032842183\\
0.43359375	0.00311856032842183\\
0.4375	0.00311856032842183\\
0.44140625	0.00311856032842183\\
0.4453125	0.00311856032842183\\
0.44921875	0.00311856032842183\\
0.453125	0.00311856032842183\\
0.45703125	0.00311856032842183\\
0.4609375	0.00311856032842183\\
0.46484375	0.00311856032842183\\
0.46875	0.00311856032842183\\
0.47265625	0.00311856032842183\\
0.4765625	0.00311856032842183\\
0.48046875	0.00311856032842183\\
0.484375	0.00311856032842183\\
0.48828125	0.00311856032842183\\
0.4921875	0.00311856032842183\\
0.49609375	0.00311856032842183\\
0.5	0.00311856032842183\\
0.50390625	0.00311856032842183\\
0.5078125	0.00311856032842183\\
0.51171875	0.00311856032842183\\
0.515625	0.00311856032842183\\
0.51953125	0.00311856032842183\\
0.5234375	0.00311856032842183\\
0.52734375	0.00311856032842183\\
0.53125	0.00311856032842183\\
0.53515625	0.00311856032842183\\
0.5390625	0.00311856032842183\\
0.54296875	0.00311856032842183\\
0.546875	0.00311856032842183\\
0.55078125	0.00311856032842183\\
0.5546875	0.00311856032842183\\
0.55859375	0.00311856032842183\\
0.5625	0.00311856032842183\\
0.56640625	0.00311856032842183\\
0.5703125	0.00311856032842183\\
0.57421875	0.00311856032842183\\
0.578125	0.00311856032842183\\
0.58203125	0.00311856032842183\\
0.5859375	0.00311856032842183\\
0.58984375	0.00311856032842183\\
0.59375	0.00311856032842183\\
0.59765625	0.00311856032842183\\
0.6015625	0.00311856032842183\\
0.60546875	0.00311856032842183\\
0.609375	0.00311856032842183\\
0.61328125	0.00311856032842183\\
0.6171875	0.00311856032842183\\
0.62109375	0.00311856032842183\\
0.625	0.00311856032842183\\
0.62890625	0.00311856032842183\\
0.6328125	0.00311856032842183\\
0.63671875	0.00311856032842183\\
0.640625	0.00311856032842183\\
0.64453125	0.00311856032842183\\
0.6484375	0.00311856032842183\\
0.65234375	0.00311856032842183\\
0.65625	0.00311856032842183\\
0.66015625	0.00311856032842183\\
0.6640625	0.00311856032842183\\
0.66796875	0.00311856032842183\\
0.671875	0.00311856032842183\\
0.67578125	0.00311856032842183\\
0.6796875	0.00311856032842183\\
0.68359375	0.00311856032842183\\
0.6875	0.00311856032842183\\
0.69140625	0.00311856032842183\\
0.6953125	0.00311856032842183\\
0.69921875	0.00311856032842183\\
0.703125	0.00311856032842183\\
0.70703125	0.00311856032842183\\
0.7109375	0.00311856032842183\\
0.71484375	0.00311856032842183\\
0.71875	0.00311856032842183\\
0.72265625	0.00311856032842183\\
0.7265625	0.00311856032842183\\
0.73046875	0.00311856032842183\\
0.734375	0.00311856032842183\\
0.73828125	0.00311856032842183\\
0.7421875	0.00311856032842183\\
0.74609375	0.00311856032842183\\
0.75	0.00311856032842183\\
0.75390625	0.00311856032842183\\
0.7578125	0.00311856032842183\\
0.76171875	0.00311856032842183\\
0.765625	0.00311856032842183\\
0.76953125	0.00311856032842183\\
0.7734375	0.00311856032842183\\
0.77734375	0.00311856032842183\\
0.78125	0.00311856032842183\\
0.78515625	0.00311856032842183\\
0.7890625	0.00311856032842183\\
0.79296875	0.00311856032842183\\
0.796875	0.00311856032842183\\
0.80078125	0.00311856032842183\\
0.8046875	0.00311856032842183\\
0.80859375	0.00311856032842183\\
0.8125	0.00311856032842183\\
0.81640625	0.00311856032842183\\
0.8203125	0.00311856032842183\\
0.82421875	0.00311856032842183\\
0.828125	0.00311856032842183\\
0.83203125	0.00311856032842183\\
0.8359375	0.00311856032842183\\
0.83984375	0.00311856032842183\\
0.84375	0.00311856032842183\\
0.84765625	0.00311856032842183\\
0.8515625	0.00311856032842183\\
0.85546875	0.00311856032842183\\
0.859375	0.00311856032842183\\
0.86328125	0.00311856032842183\\
0.8671875	0.00311856032842183\\
0.87109375	0.00311856032842183\\
0.875	0.00311856032842183\\
0.87890625	0.00311856032842183\\
0.8828125	0.00311856032842183\\
0.88671875	0.00311856032842183\\
0.890625	0.00311856032842183\\
0.89453125	0.00311856032842183\\
0.8984375	0.00311856032842183\\
0.90234375	0.00311856032842183\\
0.90625	0.00311856032842183\\
0.91015625	0.00311856032842183\\
0.9140625	0.00311856032842183\\
0.91796875	0.00311856032842183\\
0.921875	0.00311856032842183\\
0.92578125	0.00311856032842183\\
0.9296875	0.00311856032842183\\
0.93359375	0.00311856032842183\\
0.9375	0.00311856032842183\\
0.94140625	0.00311856032842183\\
0.9453125	0.00311856032842183\\
0.94921875	0.00311856032842183\\
0.953125	0.00311856032842183\\
0.95703125	0.00311856032842183\\
0.9609375	0.00311856032842183\\
0.96484375	0.00311856032842183\\
0.96875	0.00311856032842183\\
0.97265625	0.00311856032842183\\
0.9765625	0.00311856032842183\\
0.98046875	0.00311856032842183\\
0.984375	0.00311856032842183\\
0.98828125	0.00311856032842183\\
0.9921875	0.00311856032842183\\
0.99609375	0.00311856032842183\\
1	0.00311856032842183\\
};
\addlegendentry{Error bound \texttt{BT}}

\end{axis}
\end{tikzpicture}%

%% file: normH2decay_os_vs_bt.tex
%
%
\definecolor{mycolor1}{rgb}{0.00000,0.44700,0.74100}%
\definecolor{mycolor2}{rgb}{0.85000,0.32500,0.09800}%
\begin{tikzpicture}

\begin{axis}[%
width=0.951\wex,
height=\hex,
at={(0\wex,0\hex)},
scale only axis,
xmin=0,
xmax=25,
xlabel style={font=\color{white!15!black}},
xlabel={reduced order $r$},
ymode=log,
ymin=0.004,
ymax=2,
yminorticks=true,
ylabel style={font=\color{white!15!black}},
ylabel={Error bound},
axis background/.style={fill=white},
title style={font=\bfseries},
title={Error bound as function of the reduced order},
xmajorgrids,
ymajorgrids,
yminorgrids,
legend style={legend cell align=left, align=left, draw=white!15!black}
]
\addplot [color=mycolor1, line width=1.5pt,mark=*]
  table[row sep=crcr]{%
1	1.48036805922148\\
2	0.823417299511922\\
3	0.501907573897824\\
4	0.43286713744336\\
5	0.27367487286125\\
6	0.219404628838092\\
7	0.157854286510749\\
8	0.129775168103114\\
9	0.103109279439526\\
10	0.0812994082556144\\
11	0.0665023911398721\\
12	0.0512054996194745\\
13	0.0453149335521843\\
14	0.0384844884863643\\
15	0.0313704932323981\\
16	0.0238060035968209\\
17	0.0203535963415261\\
18	0.0175419430458386\\
19	0.0148429190211184\\
20	0.0122313413077736\\
21	0.0103877333585853\\
22	0.00858311584857121\\
23	0.00736758472459284\\
24	0.00627450971901284\\
25	0.00546136541788652\\
};
\addlegendentry{\texttt{OS}}

\addplot [color=mycolor2, line width=1.5pt, mark=square]
  table[row sep=crcr]{%
1	1.96124495584831\\
2	1.02909186524838\\
3	0.555386921085935\\
4	0.494413312535458\\
5	0.297114383133617\\
6	0.240808656232321\\
7	0.171208614900113\\
8	0.136619338464254\\
9	0.110637278402729\\
10	0.0880514837285449\\
11	0.0697048313870464\\
12	0.0539768095597703\\
13	0.0485190180791991\\
14	0.0412215367173302\\
15	0.0333350889696252\\
16	0.0248445000345581\\
17	0.0214123351237039\\
18	0.0182546074826567\\
19	0.0152831413556376\\
20	0.0127066864374285\\
21	0.0108622306756786\\
22	0.00892178278428066\\
23	0.00762669933766432\\
24	0.00644734252501685\\
25	0.00563180120271463\\
};
\addlegendentry{\texttt{BT}}

\end{axis}
\end{tikzpicture}%

%% file: main.bbl
\begin{thebibliography}{10}

\bibitem{typeIBT}
S.~A. Al-Baiyat and M.~Bettayeb.
\newblock {A new model reduction scheme for k--power bilinear systems}.
\newblock {\em Proceedings of the 32nd IEEE Conference on Decision and
  Control}, pages 22--27, 1993.

\bibitem{bennerdamm}
P.~Benner and T.~Damm.
\newblock {Lyapunov equations, energy functionals, and model order reduction of
  bilinear and stochastic systems.}
\newblock {\em SIAM J. Control Optim.}, 49(2):686--711, 2011.

\bibitem{redbendamm}
P.~Benner, T.~Damm, M.~Redmann, and Y.~R. Rodriguez~Cruz.
\newblock {Positive Operators and Stable Truncation.}
\newblock {\em Linear Algebra Appl}, 498:74--87, 2016.

\bibitem{bennerdammcruz}
P.~Benner, T.~Damm, and Y.~R. Rodriguez~Cruz.
\newblock {Dual pairs of generalized Lyapunov inequalities and balanced
  truncation of stochastic linear systems}.
\newblock {\em IEEE Trans. Autom. Contr.}, 62(2):782--791, 2017.

\bibitem{benner2015survey}
P.~Benner, S.~Gugercin, and K.~Willcox.
\newblock A survey of projection-based model reduction methods for parametric
  dynamical systems.
\newblock {\em SIAM Rev.}, 57(4):483--531, 2015.

\bibitem{redmannbenner}
P.~{Benner} and M.~{Redmann}.
\newblock {Model Reduction for Stochastic Systems.}
\newblock {\em {Stoch PDE: Anal Comp}}, 3(3):291--338, 2015.

\bibitem{BenS13}
P.~Benner and J.~Saak.
\newblock {Numerical solution of large and sparse continuous time algebraic
  matrix {R}iccati and {L}yapunov equations: a state of the art survey}.
\newblock {\em GAMM Mitteilungen}, 36(1):32--52, 2013.

\bibitem{morBenS10}
P.~Benner and A.~Schneider.
\newblock {Balanced Truncation Model Order Reduction for {LTI} Systems with
  many Inputs or Outputs}.
\newblock In Andr\'as Edelmayer, editor, {\em Proc. of the 19th International
  Symposium on Mathematical Theory of Networks and Systems}, pages 1971--1974,
  Budapest, Hungary, 2010.

\bibitem{berkooz1993proper}
G.~Berkooz, P.~Holmes, and J.~L. Lumley.
\newblock The proper orthogonal decomposition in the analysis of turbulent
  flows.
\newblock {\em Annual review of fluid mechanics}, 25(1):539--575, 1993.

\bibitem{breiten2010krylov}
T.~Breiten and T.~Damm.
\newblock Krylov subspace methods for model order reduction of bilinear control
  systems.
\newblock {\em Syst. Control. Lett.}, 59(8):443--450, 2010.

\bibitem{morcondon2005}
M.~Condon and R.~Ivanov.
\newblock Nonlinear systems-algebraic gramians and model reduction.
\newblock {\em COMPEL}, 24(1):202--219, 2005.

\bibitem{damm}
T.~Damm.
\newblock {\em {Rational Matrix Equations in Stochastic Control.}}
\newblock {Lecture Notes in Control and Information Sciences 297. Berlin:
  Springer}, 2004.

\bibitem{enns1984model}
D.~F. Enns.
\newblock {Model reduction with balanced realizations: An error bound and a
  frequency weighted generalization}.
\newblock In {\em The 23rd IEEE conference on decision and control}, pages
  127--132. IEEE, 1984.

\bibitem{freund2003model}
R.~W. Freund.
\newblock {Model reduction methods based on Krylov subspaces}.
\newblock {\em Acta Numerica}, 12:267--319, 2003.

\bibitem{graymesko}
W.~S. Gray and J.~Mesko.
\newblock {Energy Functions and Algebraic Gramians for Bilinear Systems}.
\newblock {\em Proceedings of the 4th IFAC Nonlinear Control Systems Design
  Symposium}, 31(17):101--106, 1998.

\bibitem{staboriginal}
R.~Z. Khasminskii.
\newblock {\em Stochastic stability of differential equations}, volume~66 of
  {\em Stochastic Modelling and Applied Probability}.
\newblock Springer, Heidelberg, second edition, 2012.

\bibitem{kunisch2001galerkin}
K.~Kunisch and S.~Volkwein.
\newblock {Galerkin proper orthogonal decomposition methods for parabolic
  problems}.
\newblock {\em Numerische Mathematik}, 90(1):117--148, 2001.

\bibitem{moo1981}
B.~Moore.
\newblock Principal component analysis in linear systems: Controllability,
  observability, and model reduction.
\newblock {\em IEEE transactions on automatic control}, 26(1):17--32, 1981.

\bibitem{pernebo1982model}
L.~Pernebo and L.~Silverman.
\newblock Model reduction via balanced state space representations.
\newblock {\em IEEE Transactions on Automatic Control}, 27(2):382--387, 1982.

\bibitem{polyuga2010model}
R.~V. Polyuga and A.~van~der Schaft.
\newblock {Model reduction of port-Hamiltonian systems as structured systems}.
\newblock In {\em Proceedings of the 19th International Symposium on
  Mathematical Theory of Networks and Systems--MTNS}, volume~5, 2010.

\bibitem{prajna2003pod}
S.~Prajna.
\newblock {POD} model reduction with stability guarantee.
\newblock In {\em 42nd IEEE International Conference on Decision and Control
  (IEEE Cat. No. 03CH37475)}, volume~5, pages 5254--5258. IEEE, 2003.

\bibitem{morPul19}
R.~Pulch.
\newblock Stability preservation in {G}alerkin-type projection-based model
  order reduction.
\newblock {\em Numer. Algebra, Control. Optim.}, 9(1):23--44, 2019.

\bibitem{redmannPhD}
M.~Redmann.
\newblock {\em {Model Order Reduction Techniques Applied to Evolution Equations
  with Lévy Noise}}.
\newblock PhD thesis, {Otto-von-Guericke-Universität Magdeburg}, 2016.

\bibitem{redmannspa2}
M.~Redmann.
\newblock {Type II singular perturbation approximation for linear systems with
  Lévy noise}.
\newblock {\em SIAM J. Control Optim.}, 56(3):2120--2158., 2018.

\bibitem{h2_bil}
M.~Redmann.
\newblock {Bilinear systems -- A new link to $\mathcal H_2$-norms, relations to
  stochastic systems and further properties}.
\newblock {\em arXiv preprint: 1910.14427v4}, 2019.

\bibitem{BTtyp2EB}
M.~Redmann and P.~Benner.
\newblock {An $H_2$-Type Error Bound for Balancing-Related Model Order
  Reduction of Linear Systems with L\'evy Noise}.
\newblock {\em Syst. Control. Lett.}, 105:1--5, 2017.

\bibitem{mliopt}
M.~Redmann and M.~A. Freitag.
\newblock {Optimization based model order reduction for stochastic systems}.
\newblock {\em Appl. Math. Comput.}, Volume 398, 2021.

\bibitem{selga2012stability}
R.~C. Selga, B.~Lohmann, and R.~Eid.
\newblock Stability preservation in projection-based model order reduction of
  large scale systems.
\newblock {\em Eur J Control}, 18(2):122--132, 2012.

\bibitem{Sim16a}
V.~Simoncini.
\newblock Computational methods for linear matrix equations.
\newblock {\em {SIAM} Rev.}, 38(3):377--441, 2016.

\bibitem{sorensen2005model}
D.~C. Sorensen and A.~C. Antoulas.
\newblock On model reduction of structured systems.
\newblock In {\em Dimension reduction of large-scale systems}, pages 117--130.
  Springer, 2005.

\bibitem{wonham}
W.~M. Wonham.
\newblock {On a Matrix Riccati Equation of Stochastic Control}.
\newblock {\em SIAM J. Control}, 6(4):681--697, 1968.

\end{thebibliography}
